\documentclass[%
  twocolumn
   , colorlinks 
]{mpi2015-cscpreprint}

\usepackage[T1]{fontenc}
\usepackage[utf8]{inputenc}
\usepackage[american]{babel}

\usepackage{amsmath,amssymb,amsfonts,mathtools}
\usepackage{stmaryrd}
\usepackage{bbold}

\usepackage[linesnumbered,ruled,vlined]{algorithm2e}

\usepackage{graphicx}
\usepackage{import}
\usepackage[export]{adjustbox} 
\usepackage{subcaption}        
\usepackage{pgfplots}
\usepgfplotslibrary{groupplots}
\pgfplotsset{compat=1.18}
\usepackage{tikz}
\usetikzlibrary{math,shapes,snakes,arrows.meta}
\usepackage{psfrag}

\usepackage{booktabs}
\usepackage{tabularx}

\usepackage{siunitx}
\usepackage{float} 

\pgfplotscreateplotcyclelist{three_col}{
	{color of colormap={0 of viridis}, mark=*, line width=0.75pt, mark repeat=6, mark phase=0, mark options={draw=black}},
	{color of colormap={300 of viridis}, mark=square*, line width=1.0pt, mark repeat=6, mark phase=1, mark options={draw=black}},
	{color of colormap={600 of viridis}, mark=triangle*, line width=1.2pt, mark repeat=6, mark phase=2, mark options={draw=black}},
}

\newcommand{\MATLAB}{\textsc{Matlab}}
\newcommand{\algoSeq}{\textit{algorithmic sequence for the analysis of counter-rotating detonation waves}\,}

\begin{document}
\title{Data Driven Modeling of Nonlinear Dynamics in a Rotating Detonation Combustor via Finite Dimensional Approximations of the Koopman Operator}

\author[$\ast$]{David Oexle}
\affil[$\ast$]{Modeling, Simulation, and Optimization of Real Processes, Technische Universität Berlin\authorcr%
\email{oexle@mpi-magdeburg.mpg.de}, \orcid{0009-0006-3901-628X}}

\author[$\dagger$]{Tobias Breiten}
\affil[$\dagger$]{Modeling, Simulation, and Optimization of Real Processes, Technische Universität Berlin\authorcr%
\email{tobias.breiten@tu-berlin.de}, \orcid{0000-0002-9815-4897}}

\author[$\dagger$]{Myles D. Bohon}
\affil[$\dagger$]{Institute of Fluid Mechanics and Technical Acoustics, Technische Universität Berlin\authorcr%
\email{m.bohon@tu-berlin.de}, \orcid{0000-0002-6724-0805}}

\shorttitle{Data Driven Modeling of Nonlinear Dynamics in a RDC}
\shortauthor{D. Oexle, T. Breiten, M.D.Bohon}
\shortdate{}


\abstract{
A Rotating Detonation Combustor (RDC) is a promising technology for increasing efficiency in propulsion and power generation applications.
The dynamics of the RDC are governed by continuously propagating detonation waves within an annular combustion chamber.
Multiple operating modes can be observed, including nonlinear interactions between counter-rotating waves and the emergence of standing wave patterns.
Koopman operator theory provides a framework to globally linearize nonlinear dynamical systems by representing their evolution in the space of observables rather than states.
In this work,
finite-dimensional approximations of the Koopman operator are constructed using variants of Dynamic Mode Decomposition (DMD) applied to high-speed video data capturing the natural flame luminosity of the detonation waves in the RDC at the Technical University (TU) Berlin.
By introducing time-delay embeddings as a dictionary of observables,
this approach overcomes the limitations of standard DMD methods,
particularly for accurate reconstruction of standing wave patterns and for capturing nonlinear interactions. 
In addition,
a technique is presented to mitigate the influence of sensor noise in the luminosity measurements.
Finally,
it is shown that the DMD-based models provide insight into the dynamics of different operating modes by decomposing the reconstructed signal into its characteristic features. 
}

\novelty{ 
	\begin{itemize}
		\item  We introduce a combination of  methods tailored to the challenges of the RDC data, which then comprise an \algoSeq.
			In particular, we propose an estimate for the number of time delay embeddings based on the physics of the operating modes under consideration in the context of the Extended Dynamic Mode Decomposition.
		\item  We demonstrate how the developed sequence can be used to separate traveling waves structure from the data set and analyze the remaining dynamics separately. 
	\end{itemize} 
}

\maketitle

\section{Introduction}\label{sec:introduction}
An engine that makes use of Pressure Gain Combustion (PGC) can increase its thermodynamic efficiency.  
In contrast to combustion systems, where heat is added to the combustion cycle at constant pressure, in a PGC device heat is added at nearly constant volume, which results in a significant pressure increase \cite{kailasanathPCG00,gutmarkPGC21}.  
A promising design that makes use of the theoretical benefits of PGC is the rotating detonation combustor (RDC), in which detonation waves are constantly traveling through an annular combustion chamber as long as new reactants are supplied \cite{bluemmnerRDC19}. 

\begin{figure}[htpb]
	\centering
	\def \svgwidth{\columnwidth}
	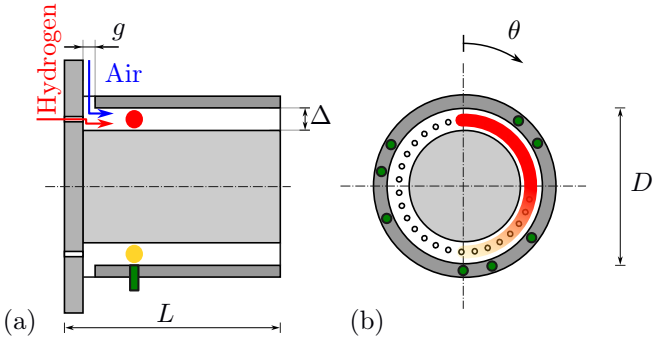
	\caption{
		Panel (a) shows a longitudinal section of the RDC and the relevant dimensions.
		Panel (b) shows a cross  section view, where the green markers indicate the position of the pressure probes.
		The fading yellow arc sketches the trajectory of a traveling detonation wave.
	}
	\label{fig:1_1_rdcSchematics}
\end{figure}

The RDC used in the studies at Technische Universität (TU) Berlin \cite{bachRDC20,bluemmnerRDC19,bohonRDE19,barnouinRDC25} is based on a design by the US Air Force Research Laboratory (AFRL) \cite{shankRDC12}, which operates on non-premixed mixtures of hydrogen and air.  
A  schematic of the RDC is shown in Figure~\ref{fig:1_1_rdcSchematics} and illustrates the  structure of the combustor. 
Further, the dimensions of the RDC are detailed in Table~\ref{tab:1_1_rdcDimensions}.
\begin{table}[htpb]
    \centering
    \caption{Dimensions of the RDC}
    \label{tab:1_1_rdcDimensions} 
    \begin{tabularx}{\columnwidth}{X c c c}
        \toprule
        Property & Symbol & Value & Unit \\
        \midrule
        Fuel hole diameter & $d$ & 0.5 & mm \\
        Number of fuel holes & $N$ & 100 & - \\
        Oxidizer slot height & $g$ & 1.6 & mm \\
        Annulus outer diameter & $D$ & 90 & mm \\
        Annulus width & $\Delta$ & 7.6 & mm \\
        Annulus Length & $L$ & 120 & mm \\
        \bottomrule
    \end{tabularx}
\end{table}

In contrast to the rather simplistic design of the 
RDC, the observed physical phenomena are complex. 
During operation, several nonidealities such as 
counter-rotating waves can occur. 
A review detailing the physical background can be 
found in \cite{ramanNonRDE2023} and \cite{anandRDC19}.

While operating the device,
different modes of operation are observed.  
Apart from the canonical case of a single spinning wave traveling around the annular combustion chamber,
various operating modes with counter-rotating wave structures \cite{barnouinRDC25} are found as well.  
Therefore,
a central challenge to gain insights into the stability and efficiency of the RDC is to understand and quantify the dynamics of nonlinear interactions such as nonlinear amplification.  
In order to model the nonlinear dynamics,
different strategies can be pursued.  
A first-order model approach, i.e.,
finding a set of differential equations via balance equations and constitutive laws,
is described in \cite{kochModelRDE20,kochMultiRDE21},
which focuses on the interaction and bifurcation behavior of co-rotating detonation waves.
The subsequent work by \cite{kochNeuralODE2021, baoCheap2RichRDE26} demonstrates how the knowledge on the governing equations can be incorporated into data-driven deep learning frameworks.
A purely data-driven approach for system identification can be found in \cite{johnsonApplicationConvolutionalNeural2021}.
Further, efforts to build surrogate models from high-fidelity simulation data consolidate model order reduction approaches such as operator 
inference \cite{farcasOpInf2023} and reduced basis methods \cite{camachoProjROM2025}.

Another perspective on nonlinear dynamical systems is given through the framework of the Koopman operator theory.  
Introduced by \cite{koopHam31},
the Koopman operator observes a nonlinear dynamical system through the lens of observable functions through which the dynamics appear linear.  
Therefore,
instead of evolving the system state nonlinearly in state space,
the evolution of functions of the state space is considered in the infinite-dimensional space of observable functions.  
Hence,
the Koopman operator theory carries the promise of a global linearization of a nonlinear dynamical system at the price of being infinite-dimensional.  
A well-established procedure for finding finite-dimensional approximations of the Koopman operator from discrete measurements of the system is the Dynamic Mode Decomposition (DMD),
introduced by \cite{schmidDMD10} and later connected to the Koopman operator theory by \cite{rowleyKoop09}.  
Furthermore,
the spectral information of the finite-dimensional operator computed by DMD allows for identifying coherent structures in the data such as oscillatory patterns.  
Therefore,
DMD is a well-suited tool for analyzing data capturing the nonlinear dynamics of the RDC,
which can be generated by imaging the exhaust flame luminosity or sampling pressure data from the combustion chamber.  
For example,
in the work by \cite{bohonDMD21},
the DMD procedure is applied to luminosity data.  
It shows how DMD can successfully be applied to identify the dynamical features of operating modes exhibiting counter-rotating waves.  
Furthermore,
it demonstrates how the presence of standing waves and sensor noise poses a challenge when building DMD models.  
Another challenge for building DMD models in the context of RDC data is traveling waves,
which are inherently difficult to describe for mode-based methods such as DMD \cite{reissSPOD18}.  
In order to address this problem,
in the work of \cite{mendibleRDE21},
mainly observing operating modes with multiple co-rotating waves,
the dataset is transformed into a moving frame of reference (FR) prior to the DMD analysis.  
The transformation from the inertial into the moving frame of reference is identified via a learning algorithm as proposed in \cite{mendibleROMwave20}.  
Another shortcoming of the standard DMD method is that it cannot truly depict nonlinear phenomena.  
To this end,
\cite{williamsEDMD15} introduced the Extended Dynamic Mode Decomposition (EDMD),
which provides a more general framework for finding finite-dimensional approximations of the Koopman operator by using a sequence of snapshots and a dictionary of functions.

The goal of this work is to develop  a workflow for the analysis of the nonlinear dynamics of counter rotating detonation waves in a RDC,
based on finite-dimensional approximations of the Koopman operator.

In particular, 
we introduce a combination of methods tailored to the challenges of the RDC data,
which then comprise an \algoSeq.
Further, we propose an estimate for the number of 
time delay embeddings based on the physics of the operating modes under consideration in the context 
of the Extended Dynamic Mode Decomposition.
We demonstrate how the developed sequence can be used to recover precise estimates of the counter 
rotating traveling wave structures from the DMD analysis in the moving frame of reference.
Further,
separating the traveling wave structures 
from the data set allows us to analyze the remaining dynamics separately.
Hence, we offer a workflow that enables 
RDE practitioners to analyze the modal structure of 
counter-rotating detonation waves within the DMD framework.
Therefore, the methodological novelty of this work lies more in the problem specific composition of existing DMD variants rather than in the advances of the DMD method itself.

The remainder of this work is structured as follows.  
In section \ref{sec:setup},
the setup of the RDC at TU Berlin is presented and it is detailed how image data of the flame luminosity and pressure data are collected.  
Furthermore,
a data post-processing routine is described in order to reduce the dimensionality of the image data describing the luminosity.  
Subsequently,
the different operating modes and their distinct dynamical features are then presented.  
In  section \ref{sec:methods}, 
the methodological concepts that comprise the sequence are introduced.  
Firstly,
the Koopman operator theory is introduced as a theoretical framework for analyzing nonlinear dynamics.  
Subsequently,
it is shown how finite-dimensional approximations of the Koopman operator can be obtained using the EDMD procedure from a sequence of data and a set of dictionary functions. 
Additionally,
the Exact Dynamic Mode Decomposition (ExactDMD) \cite{tuDMD14} is used as an algorithmic variant for computing the desired finite-dimensional approximation of the Koopman operator.  
In order to remove the bias introduced by the presence of sensor noise in the data, the Optimized Dynamic Mode Decomposition (OptDMD) \cite{ashkamNoiseOptDMD17} algorithm is applied as a post-processing step for the results of the ExactDMD.  
Furthermore, a criterion that separates the noise from the actual dynamics \cite{gavishSigTresh14} is considered.  
Finally,
a procedure is introduced to change the frame of reference (FR) according to the wave speed of the dynamical features of the operating mode.
In section \ref{sec:results}, it is demonstrated how the standard DMD (exactDMD) method fails when it is applied to the RDC luminosity data. 
Subsequently,
the algorithmic steps of the \algoSeq are presented 
and applied to the luminosity data of the introduced operating modes of the RDC.
Finally,
section \ref{sec:discussion} concludes the work with a brief summary of the results and an outlook on future work.

\section{Setup}\label{sec:setup}
\subsection{Experimental setup}\label{sec:experimental-setup}
\begin{figure*}[htpb]
	\centering
	\def \svgwidth{0.8\textwidth}
	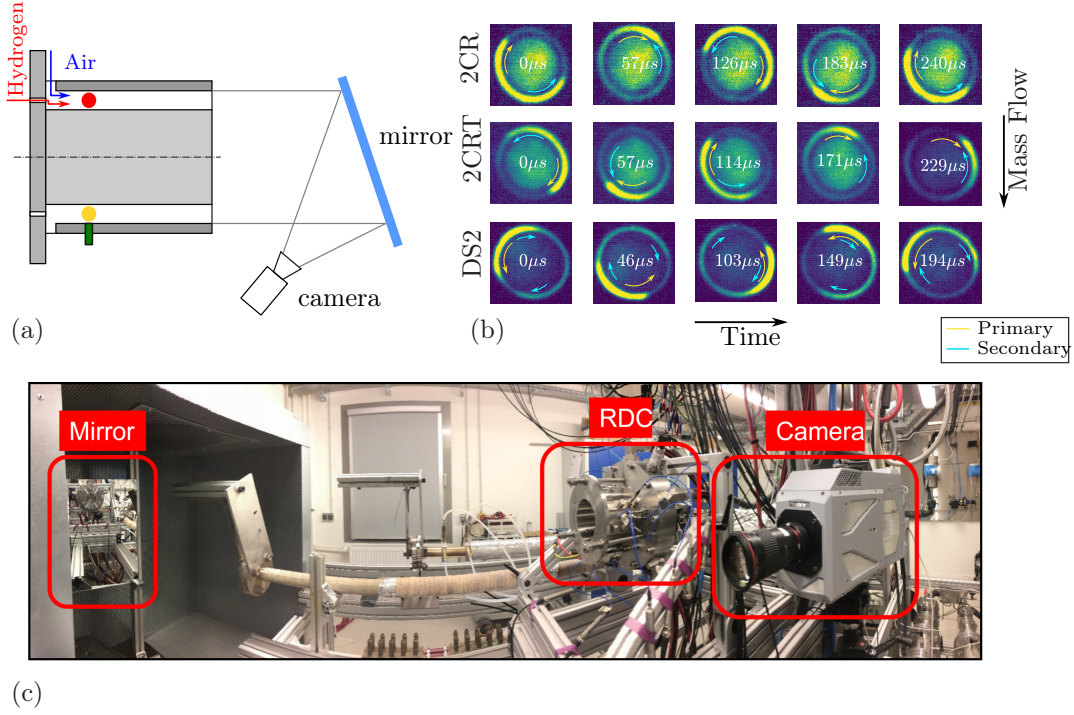  
	\caption{The left panel (a) of the figure shows the experimental setup for imaging the aft end of the RDC.  
		The right panel (b) shows a series of snapshots of the three considered operating modes generated by increasing the mass flow of the reactants in the RDC.  
		Each row corresponds to one operating mode,
	illustrating its temporal evolution over one full period of wave propagation. Panel (c) shows a photograph of the experimental setup.}   
	\label{fig:2_1_experimentalSetup}
\end{figure*}

In order to gather data on the dynamical processes of detonation waves traveling and colliding within the combustion chamber,
two data acquisition methods are introduced.  
The first method focuses on collecting pressure data via piezoelectric pressure sensors (PCB).  
Figure~\ref{fig:1_1_rdcSchematics}(b) indicates the azimuthal positioning of eight PCBs distributed along the annular combustion chamber at a fixed height close to the base end of the combustion chamber.
The sensors resolve the pressure data over time at a data acquisition rate of $500\,\si{\kilo\hertz}$.  
Complementing the pressure data,
image data of the flame luminosity is also gathered. 
To this end,
a high-speed camera (Photron SA-Z) is used to image the aft end of the annulus.  
To protect the camera setup from the harsh RDC exhaust,
a visible wavelength mirror is used for imaging the flame luminosity.  
However,
during operation,
the mirror deflects,
which introduces distortions in the data that must be accounted for when working with the luminosity recordings.  
The images are recorded at a frame rate of $f_s = 87{,}500$ frames per second with an exposure time of \SI{8.75}{\micro \second} and an inter-frame time of \SI{11.4}{\micro \second}.
Figure~\ref{fig:2_1_experimentalSetup}(a) shows a sketch of the experimental setup used for imaging the detonation waves in the combustion chamber.
Across a range of reactant mass flows and equivalence ratios,
different stable operating modes are observed in the RDC \cite{barnouinRDC25}.
The parameters for the experimental runs are summarized in Table \ref{tab:experimental_parameters}. 
In particular,
next to the canonical case of a single spinning wave,
a variety of operating modes exhibiting counter-rotating wave patterns are observed,
which are detailed in  Figure~\ref{fig:2_1_experimentalSetup}(b) over the course of one period.
The top row in Figure~\ref{fig:2_1_experimentalSetup}(b),
labeled (2CR),
shows the operating mode with two counter-rotating waves of similar strength and  equal speed.
This operating mode features a fixed azimuthal interaction point and generally occurs at low mass flow rates or high outlet blockage ratios.
In particular,
this mode results in a standing wave pattern.
Increasing the mass flow leads to one of the counter-rotating waves increasing in velocity and amplitude relative to the other.
This results in an operating mode characterized by a dominant primary wave and a slower secondary wave,
with a precessing interaction point.
This mode is shown in the second row of  Figure~\ref{fig:2_1_experimentalSetup}(b),
and is labeled (2CRT).
Further increasing the mass flow leads to a continued amplification and acceleration of the primary wave.
Meanwhile,
the secondary wave tends to split into a set of two or more weaker and slower counter-rotating waves.
This mode is illustrated in the third row of Figure~\ref{fig:2_1_experimentalSetup}(b) and is labeled (DS2).

\begin{table}[htpb]
	\caption{Overview of experimental parameters across datasets for different RDC operating modes}
	\label{tab:experimental_parameters}
	\begin{tabularx}{\columnwidth}{l c c c c c c}
			\toprule
			Property & Symbol & 2CR & 2CRT & DS2 &  Unit \\
			\midrule 
			Mass flow & $\dot{m}$ & 123.8 & 147.5 & 199.1 &  \si{g/s} \\
			Equiv. ratio  \footnote{The equivalence ratio $\phi$ describes the ratio of the fuel-to-oxidizer ratio to the stoichiometric fuel-to-oxidizer ratio.} & $\phi$  &
			{0.90}&
			{0.89}&
			{1.00}&
			-- \\
			Runtime & $t$ & 57 & 57 & 57 & \si{ms} \\
			Time steps & $m$ & 5000 & 5000  &5000 & - \\
			\bottomrule
		\end{tabularx}
\end{table}

\subsection{Data pipeline}\label{sec:data-pipeline}
\begin{figure}[htpb]
	\centering
	\def \svgwidth{\columnwidth}
	\small{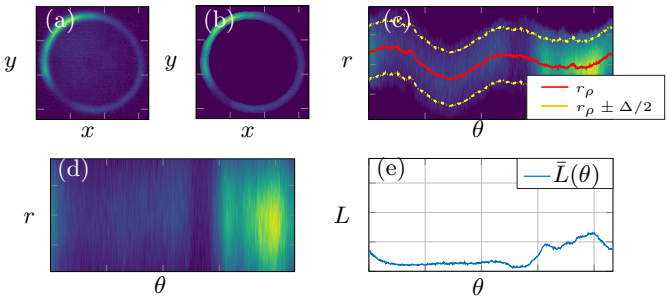}
	\caption{Panel (a)  shows the unprocessed image. Panel (b) shows the  masked image obtained by a luminosity threshold. Panel (c) shows the image in polar coordinates, the mean center of luminosity $r_\rho$ and the bounds $r_\rho \pm \frac{\Delta}{2}$. Panel (d) shows the centered data. Panel (e) shows the radially averaged luminosity  $\bar{L}(\theta)$.}
	\label{fig:DataPipeline}
\end{figure}

In a single recording session,
up to $5000$ images are generated,
resolving the natural flame luminosity over a grid of up to $329\times329$ pixels.  
This results in a very high dimensional representation of the system state at a given point in time.  
However,
as illustrated in Figure~\ref{fig:2_1_experimentalSetup},
which shows the different operating modes,
the dynamics of the traveling detonation waves are confined to a subset of each image snapshot that corresponds to the annular cross section of width $\Delta$.  
In addition,
	it is noted that the data acquisition via the high-speed video camera effectively integrates the luminosity along the length of the combustion chamber.
	This limits the capability of the imaging setup to resolve small-scale features along the radial dimension.
	Therefore,
	it is assumed that the luminosity does not vary significantly across the radial direction within the annulus.
Furthermore,
the image snapshots show that the detonation waves and their trails do not follow a perfectly circular path.  
As alluded to earlier,
this is due to the exhaust flow deflecting the mirror used for imaging the aft end of the RDC.  
In order to account for this distortion,
and to exploit the fact that the dynamics are localized to a subset of the image,
a sequence of post-processing steps for the luminosity data is developed in a data pipeline.  
The individual steps of this pipeline are detailed in Figure~\ref{fig:DataPipeline} for an exemplary snapshot.

In conclusion,
the data post-processing routine,
as presented for one example in the data pipeline,
reduces the data dimensionality from image snapshots of size $329 \times 329$ pixels to snapshot vectors $\boldsymbol{z}_k \in \mathbb{R}^n$ where $n=819$. 
Furthermore,
the data is centered according to the actual width of the annular gap.  

\subsection{Dynamical feature of operating modes}\label{sec:analysis-operating-mode}
\begin{figure*}[htpb]
	\centering
	\def \svgwidth{0.85\textwidth}
	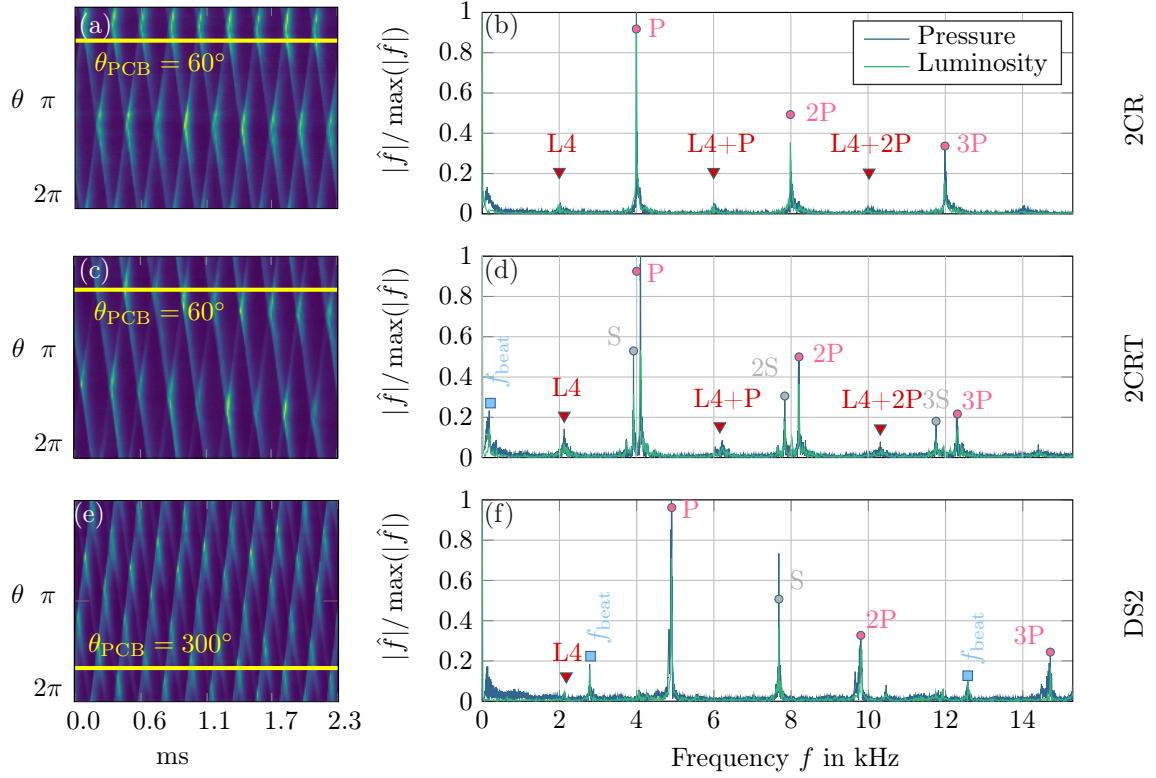
	\caption{Panel (a), (c) and (e) show the $\theta$- $t$ diagrams of the different operating modes, 
		where the yellow lines indicate the azimuthal positions  $\theta_\mathrm{PCB}.$  
	Panel (b), (d) and (f) compare the FFT of the pressure signal with the FFT of the luminosity signal. }
	\label{fig:2_3_dynamical_features}
\end{figure*}

With the dimensionality reduced and centered data vectors at hand,
they can now be stacked to form a $\theta $-$t$ diagram,
which allows for better identification of the dynamical features associated with each operating mode.
To quantify and distinguish the dynamical features,
a Fast Fourier Transform (FFT) is applied to the processed luminosity signal at the azimuthal position $\theta_{\mathrm{PCB}}$,
matching that of the corresponding pressure probe.  
Subsequently,
the peaks in the resulting frequency spectrum are assigned to the dynamical features of the operating mode,
including primary waves (P),
secondary waves (S),
longitudinal waves (L4),
and their interactions and harmonics. 
Figure~\ref{fig:2_3_dynamical_features} (a), (c) and (e) show a subset of the data set representing the operating mode in  a  $\theta$-$t$ diagram accordingly.
Furthermore,
a yellow line indicates the azimuthal position $\theta_\mathrm{PCB}$ of the luminosity probe.
Figure~\ref{fig:2_3_dynamical_features} (b), (d) and (f) show the normalized spectrum of the FFT of the pressure and luminosity signals sampled at the  azimuthal position.
The individual peaks are labeled according to the assigned feature dynamical feature.          
For the 2CR case,
the two counter rotating waves are identified by one peak at around \SI{4}{kHz},
labeled (P), and its integer multiples ($i \times \mathrm{P}$).  
Since for this operating mode the two counter rotating waves travel at identical speeds,
they cannot be distinguished from one another in the frequency spectrum.  
The presence of a weak longitudinal wave is indicated by the peak at $2\si{kHz}$, labeled (L4),
and by its interactions with the primary wave and its multiples ($ i \times \mathrm{P} + \mathrm{L4}$).
In the 2CRT case,
the counter-rotating waves are now traveling at different speeds along the combustion chamber.  
As a result,
the interaction point slowly drifts over time,
as can be observed in the $\theta$-$t$ diagram.
The faster primary wave is identified by the peak at 4.1 \si{kHz} labeled (P) and its integer multiples ($i \times \mathrm{P}$).  
The slower counter-rotating wave,
now distinguishable from the primary wave,
is labeled (S) and corresponds to the peak at 3.9 \si{kHz} and its integer multiples ($i \times \mathrm{S}$).  
The difference in wave speeds results in a beating frequency $f_\mathrm{beat} = f_\mathrm{P} - f_\mathrm{S}$,
indicated by the peak at 0.2 \si{kHz}.
Similar to the 2CR operating mode,
the presence of a longitudinal wave is indicated by the peak at around 2 \si{kHz}, labeled L4.  
Furthermore,
interactions of the longitudinal wave with the primary $i \times \mathrm{P} + \mathrm{L4}$ are present in the spectrum.
For the DS2 operating mode the primary wave is indicated by the peak at 4.9\,\si{kHz}, labeled (P), and its integer multiples,
and appears faster compared to the 2CR and 2CRT operating modes.  
The set of two counter-rotating waves is indicated in the frequency spectrum by the peak at $7.7\,\si{kHz}$, labeled (S),
and its integer multiples ($i \times \mathrm{S}$).  
Since the counter-rotating dynamics in this case consist of a  set of two  waves,
the associated spectral component appears at a high  frequency compared to the other operating modes.  
The peak at $2.8\,\si{kHz}$ indicates the beating frequency $f_\mathrm{beat} =  f_\mathrm{S} - f_\mathrm{P}$.    
Furthermore,
similar to the previous cases,
a longitudinal mode is present,
indicated by the peak at $2.1\si{kHz}$,
labeled (L4).  

\section{Methods}\label{sec:methods}
In this section,
we collect the mathematical framework required to build the \algoSeq,
which is composed of a series of DMD-like methods and data processing steps.
We aim to tailor the individual steps of the algorithmic sequence  such that the various challenges exhibited by the RDC data (e.g.  nonlinearity, sensor noise, and traveling wave physics) are addressed accordingly.

\subsection{Nonlinearity}
\subsubsection{Koopman operator theory} \label{sec:koop-operator-theory}
The dynamics of the RDC are characterized by nonlinear phenomena,
such as nonlinear wave amplification occurring at the interaction point.
One approach for analyzing such nonlinearities within the context of dynamical systems is provided by the Koopman operator framework,
in which the evolution of a nonlinear system is observed through the lens of an observable function,
under which the dynamics evolve linearly.
In the setting of the RDC data,
one only has access to discrete measurements of the system sampled directly from the state trajectory at time step $k$, i.e., $\boldsymbol{z}_{k}= \boldsymbol{z}\left( t_k \right)  =  \boldsymbol{z} \left( k\Delta t \right)$ \cite{bruntonKoop21}.
The discrete measurements $\boldsymbol{z}_{k}$ of the system are governed by the time discrete dynamical system $\boldsymbol{z}_{k+1} = \boldsymbol{F}(\boldsymbol{z}_{k})$, 
where $\boldsymbol{F}: \mathbb{R}^n \rightarrow \mathbb{R}^n$ denotes the nonlinear flow map.  
Instead of acting on the state $\boldsymbol{z}$,
the Koopman operator acts linearly \cite{bruntonBook21}  on a vector valued observable function  $\boldsymbol{g}$ as

\begin{equation}
	\mathcal{K} \boldsymbol{g}(z)=\boldsymbol{g}\left(\boldsymbol{F}\left(z_k\right)\right)=\boldsymbol{F}\left(z_k\right)=\sum_{j=1}^N \boldsymbol{v}_j \lambda_j \varphi_j\left(z_k\right).
\end{equation} 
The derivation of the decomposition is detailed in  Appendix~\ref{sec:koop-operator-theory-appendix}. 
Hence,
by lifting the dynamics into the space of observable functions,
one trades a nonlinear finite-dimensional system for a possibly infinite-dimensional yet linear system.  
The tuples $\left\{ \left( \lambda_j, \varphi_j, \boldsymbol{v}_j \right) \right\}_{j=1}^{N}$ composed of the Koopman eigenvalues,
Koopman eigenfunctions and the Koopman modes give rise to a linear representation of the evolution of the full state governed by nonlinear dynamics \cite{williamsEDMD15}. 

Although the Koopman operator theory provides a powerful framework that carries the promise of global linearization of nonlinear dynamical systems,
there remains the main challenge that the Koopman operator is infinite-dimensional and therefore computationally intractable.
Furthermore, 
to this point it remains unclear how to actually compute the desired tuples,
that are required for the desired linear representation of the state evolution. 

\subsubsection{Extended dynamic mode decomposition} \label{sec:extended-dmd}
In order to circumvent the problematic infinite-dimensional character of the Koopman operator, finding finite-dimensional approximations and the associated tuples is of particular interest.

The EDMD,
introduced by \cite{williamsEDMD15},
aims to generate those approximation of the Koopman operator,
and the associated tuples $\left\{\left(\lambda_j, \varphi_j, \boldsymbol{v}_j\right)\right\}_{j=1}^{d}$ from a given sequence of snapshot data $\left\{\boldsymbol{z}_0, \ldots, \boldsymbol{z}_m\right\}$ and a chosen set of dictionary functions  
\begin{equation}\label{eq:dict}
	\mathcal{D} = \left\{\psi_1, \psi_2, \ldots, \psi_d\right\},
\end{equation}
where $\psi_i \in \mathcal{F}$ and $\mathrm{span}(\mathcal{D}) = \mathcal{F}_D \subset \mathcal{F}$.

First the elements contained in the dictionary are sorted into a vector valued function $\boldsymbol{\psi}: \mathcal{M} \rightarrow \mathbb{R}^{1 \times d}$ such that
\begin{equation}
	\boldsymbol{\psi}(\boldsymbol{z}) = \left[\psi_1(\boldsymbol{z}), \psi_2(\boldsymbol{z}), \ldots, \psi_{d}(\boldsymbol{z})\right].
\end{equation}
Therefore,
by definition,
any function in the span of the dictionary,
i.e.,
$\phi \in \mathcal{F}_{\mathcal{D}}$,
can be expanded as
\begin{equation}
	\phi(\boldsymbol{z}) = \sum_{j=1}^{d} a_j \psi_j(\boldsymbol{z}) = \boldsymbol{\psi}(\boldsymbol{z}) \cdot \boldsymbol{a}, 
\end{equation}
where $\boldsymbol{a} \in  \mathbb{R}^d$ is the vector of basis coefficients. 
The action of the Koopman operator on the function $\phi$ is then given as
\begin{equation}\label{eq:koop-finiteaction}
	\mathcal{K} \phi(\boldsymbol{z}) = \boldsymbol{\psi}(\boldsymbol{F}(\boldsymbol{z})) \cdot \boldsymbol{a} = \boldsymbol{\psi}(\boldsymbol{z}) \cdot (\boldsymbol{K} \boldsymbol{a}) + r(\boldsymbol{z}),
\end{equation}
where the residual term $r \in \mathcal{F}$ compensates for the fact  that $\mathcal{F}_D$ typically does not span a Koopman invariant subspace.
This highlights how the prior choice of the dictionary determines the quality of the approximation of the Koopman operator.   
Furthermore,
\eqref{eq:koop-finiteaction} shows how in the finite-dimensional setting the action of the Koopman operator on a function $\phi$ is given as the action of a matrix $\boldsymbol{K} \in \mathbb{R}^{d \times d}$ on a vector of basis coefficients $\boldsymbol{a}$ used to express the function in the chosen dictionary.

Finally,
the finite-dimensional approximation of the Koopman operator,
i.e.,
$\boldsymbol{K}$,
can be identified as the matrix minimizing the residual $r(\boldsymbol{z})$ over all data points $\boldsymbol{z}_k$.
The solution of this optimization problem is the minimizer of the cost function
\begin{equation}\label{min-res}
	J = \frac{1}{2} \sum_{k=0}^{m-1} \left|r\left(\boldsymbol{z}_k\right)\right|^2 = \frac{1}{2} \sum_{k=0}^{m-1} \left|\left(\boldsymbol{\psi}\left(\boldsymbol{z}_{k+1}\right) -\boldsymbol{\psi} \left(\boldsymbol{z}_k\right) \boldsymbol{K}\right) \boldsymbol{a}\right|^2.
\end{equation}
By introducing the data matrices
\begin{equation}
	\begin{array}{ c c c}
		\boldsymbol{\Psi}_X =
		\begin{bmatrix} 
			\boldsymbol{\psi}\left(\boldsymbol{z}_0\right) \\
			\boldsymbol{\psi}\left(\boldsymbol{z}_1\right) \\
			\vdots \\
			\boldsymbol{\psi}\left(\boldsymbol{z}_{m-1}\right)
		\end{bmatrix} 
	& 
	\mathrm{and}

	&
	\boldsymbol{\Psi}_{Y} =
	\begin{bmatrix} 
		\boldsymbol{\psi}\left(\boldsymbol{z}_1\right) \\
		\boldsymbol{\psi}\left(\boldsymbol{z}_2\right) \\
		\vdots \\
		\boldsymbol{\psi}\left(\boldsymbol{z}_m\right)
	\end{bmatrix} 
	\end{array},
\end{equation}
where $\boldsymbol{\Psi}_X$ and $\boldsymbol{\Psi}_Y \in \mathbb{R}^{m \times d}$,
the minimization problem \eqref{min-res} can be written as a least squares problem
\begin{equation}\label{eq:koop-lsq}
	\boldsymbol{K} = \underset{\boldsymbol{\hat{K}} \in \mathbb{R}^{d \times d}}{\arg\min} \| \boldsymbol{\Psi}_Y - \boldsymbol{\Psi}_X \boldsymbol{\hat{K}} \|^2_{F}, 
\end{equation} 
where $\| \cdot \|_F$ denotes the Frobenius norm. 
The solution of \eqref{eq:koop-lsq} is then given via the Moore–Penrose pseudoinverse as follows
\begin{equation}\label{eq:lsq-problem}
	\boldsymbol{K} = \boldsymbol{\Psi}_X^{+} \boldsymbol{\Psi}_Y.
\end{equation}
As shown in \cite{williamsEDMD15} the Koopman eigenfunctions are given by $\varphi_{j}(\boldsymbol{z}) = \boldsymbol{\psi}(\boldsymbol{z}) \cdot \boldsymbol{\xi}_{j}$,
where $\boldsymbol{\xi}$ are the right eigenvectors of the matrix $\boldsymbol{K}$.
Finally,
the Koopman modes are given as $\boldsymbol{V} = (\boldsymbol{W}^H \boldsymbol{B})^T$,
where $\boldsymbol{W}$ is the matrix compromised of the left eigenvectors of $\boldsymbol{K}$,
which satisfy $\boldsymbol{w}_j^H \boldsymbol{\xi}_j = 1$ and the  matrix $\boldsymbol{B} \in \mathbb{R}^{d \times n}$,
which contains the weights to build the full state observable from the dictionary.  

\subsubsection{Time delay embedding}\label{sec:time-delay-embedd}
The central question that still remains open is the choice of dictionary $\mathcal{D}$ which is required for the EDMD procedure.
Popular approaches include Radial Basis functions, Hermite polynomials, or Fourier modes \cite{williamsKerDMD15}.
However, due to the periodic nature of the RDC data, time delay embedding of the state offers a suitable choice \cite{colbrookHHDMD24} for building a dictionary.  
To this end, from the given data sequence the dictionary is constructed as
\begin{equation}
	\begin{split}
		\mathcal{D}_{\mathrm{TDE}}& = \left\{\psi\left( \boldsymbol{z}_k \right), \mathcal{K} \psi\left( \boldsymbol{z}_k \right), \ldots, \mathcal{K}^s \psi\left( \boldsymbol{z}_k \right)  \right\} = \ldots \\
		\ldots &=\left\{\boldsymbol{\psi}\left(\boldsymbol{z}_k\right), \boldsymbol{\psi}\left(\boldsymbol{z}_{k+1}\right), \ldots, \boldsymbol{\psi}\left(\boldsymbol{z}_{k+s}\right)  \right\},
	\end{split}
\end{equation}
which, in the case of periodic data, becomes Koopman invariant, i.e.,
\begin{equation}
	\mathrm{span} \left( \mathcal{K}\mathcal{D}_{\mathrm{TDE}} \right) =  \mathrm{span} \left( \mathcal{D}_{\mathrm{TDE}} \right).
\end{equation}
Furthermore, the dictionary is specified by choosing $\boldsymbol{\psi}\left(  \boldsymbol{z} \right) = \boldsymbol{z}$ which yields the extended and transposed data matrices
\begin{equation}
	\boldsymbol{X} = \begin{bmatrix}
		\boldsymbol{z}_{0} & \boldsymbol{z}_{1} &  \dots & \boldsymbol{z}_{m-1}  \\
		\boldsymbol{z}_{1} & \boldsymbol{z}_{2} & \dots  & \boldsymbol{z}_{m} \\
		\vdots  & \vdots & & \vdots \\
		\boldsymbol{z}_{s-1}  & \boldsymbol{z}_{s} & \dots & \boldsymbol{z}_{m+s-2}  
	\end{bmatrix} \in \mathbb{R}^{ns \times m}
\end{equation}
and
\begin{equation}
	\boldsymbol{Y} = \begin{bmatrix}
		\boldsymbol{z}_{1} & \boldsymbol{z}_{2} &  \dots & \boldsymbol{z}_{m}  \\
		\boldsymbol{z}_{2} & \boldsymbol{z}_{3} & \dots  & \boldsymbol{z}_{m+1} \\
		\vdots  & \vdots & & \vdots \\
		\boldsymbol{z}_{s}  & \boldsymbol{z}_{s+1} & \dots & \boldsymbol{z}_{m+s-1}  
	\end{bmatrix} \in \mathbb{R}^{ns \times m}.
\end{equation}

The number of time delay embeddings of the state under which the span of the dictionary becomes Koopman invariant depends on the system under consideration.
For systems that evolve in a periodic manner such as the RDC data,
the amount of time delay embeddings is closely connected to the periodicity of the dynamics described in the data and is determined via 
\begin{equation}\label{eq:no-time-delay}
	s=\left\lfloor\frac{T}{\Delta t}\right\rfloor=\left\lfloor\frac{1}{f \Delta t}\right\rfloor,
\end{equation}
where $\Delta t$ is the time in between snapshots and $T$ is the time spanned to cover a full period of the data and can be determined from a Fast Fourier Transform (FFT).

\subsubsection{Exact dynamic mode decomposition}
Subsequently,
the Exact Dynamic Mode Decomposition (ExactDMD)  introduced by \cite{tuDMD14} is  used  as an algorithmic variant to compute a finite-dimensional representation of the Koopman operator on a subspace,
in order to leverage the low-rank structure in the data. 
The ExactDMD outputs the DMD eigenvalue matrix $\boldsymbol{\Lambda}  \in \mathbb{C}^{r \times r}$, the DMD modes
\begin{equation}
	\boldsymbol{\Phi}= \begin{bmatrix} \boldsymbol{\phi}_1 &  \ldots &  \boldsymbol{\phi}_r \end{bmatrix}  \in \mathbb{C}^{ns \times r}
\end{equation}
and the weight coefficients  $\boldsymbol{b} \in \mathbb{C}^r$. 
The system full state can then be reconstructed via 
\begin{equation}
	\boldsymbol{z}_k =\sum_{j=1}^r\left(\boldsymbol{B}^T \boldsymbol{\phi}_j\right) \lambda_j^k b_j,
\end{equation}
which echoes  \eqref{eq:fullstate-evolution-eig}.
The Koopman modes are approximated  by $\boldsymbol{v}_j = \boldsymbol{B}^\top \boldsymbol{\phi}_j$,
the DMD eigenvalues $\lambda_j$ are approximations of the Koopman eigenvalues. 
The $j$-th weight coefficients $b_j= \varphi_j\left(\boldsymbol{z}_0 \right)$ is the $j$-th Koopman eigenfunction evaluated  at the initial state $\boldsymbol{z}_0$.

\subsection{Optimized dynamic mode decomposition}\label{sec:optimized-DMD}
One of the major drawbacks of the ExactDMD algorithm is its sensitivity to sensor noise in the data \cite{bagheriNoiseDMD14}, which leads to a biased distribution of the identified DMD eigenvalues. In particular, the real parts of the eigenvalues tend to be overestimated, resulting in a rapid decay of the reconstructed signal. As a consequence, the ExactDMD procedure often fails to even reconstruct the training data accurately.
Several variants have been introduced to mitigate this issue, including total least squares DMD \cite{hematiNoiseDMD17}, forward-backward DMD \cite{dawsonNoiseDMD16}, and OptDMD \cite{ashkamNoiseOptDMD17,sashidarBOPDMD22}. Among these, the OptDMD has demonstrated significant improvement in robustness.
The central idea behind OptDMD is to reframe standard DMD as an exponential fitting problem, i.e,
\begin{equation}
	\underset{\boldsymbol{\omega}, \boldsymbol{\Phi}_{\boldsymbol{b}}}{\operatorname{argmin}}\left\|\boldsymbol{X}-\boldsymbol{\Phi}_{\boldsymbol{b}} \boldsymbol{T}(\boldsymbol{\omega})\right\|_F,
\end{equation}
where the matrix of weighted modes is defined as
\begin{equation}
	\boldsymbol{\Phi}_{\boldsymbol{b}} = \boldsymbol{\Phi} \, \operatorname{diag}(\boldsymbol{b}) = 
	\begin{bmatrix}
		\mid & & \mid \\
		\boldsymbol{\phi}_1 & \cdots & \boldsymbol{\phi}_r \\
		\mid & & \mid
	\end{bmatrix}
	\begin{bmatrix}
		b_1 & & \\
		    & \ddots & \\
		    & & b_r
	\end{bmatrix},
\end{equation}
and the temporal evolution matrix is given by
\begin{equation}
	\boldsymbol{T}(\boldsymbol{\omega}) =
	\begin{bmatrix}
		e^{\omega_1 t_0} & \cdots & e^{\omega_1 t_{m-1}} \\
		\vdots & \ddots & \vdots \\
		e^{\omega_r t_0} & \cdots & e^{\omega_r t_{m-1}}
	\end{bmatrix}.
\end{equation}
In \cite{ashkamNoiseOptDMD17}, this nonlinear least squares problem is solved via the variable projection method \cite{rustVarPro13}, which separates the linear and nonlinear components during optimization.
A potential drawback of the iterative procedure is the lack of guaranteed convergence. However, the performance can be significantly improved by supplying an initial guess.
To this end, the DMD eigenvalues computed from ExactDMD algorithm are used as an initial estimate and then iteratively refined using the OptDMD method to reduce the noise induced bias.
In this work, the \MATLAB\, implementation provided by \cite{ashkamNoiseOptDMD17} was used to solve the optimization problem.  

\subsection{Optimal threshold criterion}\label{sec:optimal-treshold}
Both the ExactDMD and OptDMD algorithms require a choice for the parameter $r$ indicating where to truncate the singular value decomposition
\begin{equation}
	\boldsymbol{X} \approx \boldsymbol{U}_r \boldsymbol{\Sigma}_r  \boldsymbol{V}_r^H
\end{equation}
of the extended data matrix.
Choosing a suitable rank $r$ improves the conditioning of the singular value decomposition and serves as a means to discriminate between true dynamical content and modes dominated by noise,
which are typically low in energy and therefore associated with small singular values.
Determining the truncation rank for a given dataset is often done heuristically,
for instance by examining the singular value decay and locating the bend in the spectrum \cite{bruntonBook21},
or by retaining a specified percentage based on the energy captured in the modes.

A more rigorous method is proposed in \cite{gavishSigTresh14},
where an optimal threshold $\tau$ is computed such that all singular values $\sigma_i < \tau$ are interpreted as noise and excluded.
For a rectangular matrix $ \boldsymbol{X} \in \mathbb{R}^{n \times m}$ with unknown noise level,
the threshold is given by
\begin{equation}
	\tau = \omega(\beta) \sigma_{\mathrm{med}},
\end{equation}
where $\sigma_{\mathrm{med}}$ is the median of the singular values and $\omega(\beta)$ is a  function depending only on the aspect ratio $\beta = \frac{n}{m}$ of the data matrix.
In this work,
the \MATLAB\, implementation provided by \cite{gavishSigTresh14} was used to determine the optimal truncation threshold.

\subsection{Change of observer}\label{sec:change-observer}
The transport-dominated physics of traveling waves poses a known challenge for constructing mode-based models \cite{mendibleROMwave20,reissSPOD18}. In the case of the RDC data, it is advantageous to choose a moving FR, as all traveling wave structures will appear stationary in it.   
In order to shift an arbitrary state vector at time step $j$
\begin{equation}
	\boldsymbol{z}_{j,\mathrm{in}}= \begin{bmatrix}
		z_{1} & 
		z_{2} &  
		\ldots 
		z_{n}
	\end{bmatrix}^\top
\end{equation}
into the FR of the primary wave, $h$ circular shifts are applied to the components of the state vector such that
\begin{equation}
	\boldsymbol{z}_{j,\mathrm{mv}}= \begin{bmatrix}
		z_{n-(h-1)}   &
		\ldots &
		z_{n} &
		z_{1} &
		\ldots &
		z_{n-h}
	\end{bmatrix}^\top.
\end{equation}
The number of required shifts $h$ is determined as 
\begin{equation}\label{eq:amount-shifts}
	h = \left\lfloor \frac{\Delta \theta_{c}}{\Delta \theta_n} \right\rfloor
\end{equation}
where $\Delta \theta_{c} = c \Delta t$ is the angular distance the wave travels during one time step $\Delta t$ at speed $c$, and $\Delta \theta_n = \theta_{j+1} - \theta_j$ is the spatial resolution of the azimuthal angle, as discussed in Section\,\ref{sec:data-pipeline}.
In general,
$\Delta \theta_{c}$ and $\Delta \theta_n$ are not exact multiples of one another.
Therefore, the rounding applied in \eqref{eq:amount-shifts} leads to small oscillations in the moving FR, having a similar effect on the DMD eigenvalues as sensor noise.  
Moreover, the wave speed $c$ is not known in advance and must be determined from the data.
This is accomplished by determining the frequency $f$ of the traveling waves from the FFT of a pressure or luminosity signal, as discussed in Section~\ref{sec:analysis-operating-mode}. 
The wave speed is then determined as
\begin{equation}
	c = 2 \pi f,
\end{equation}
where $f$ can be either $f_\mathrm{CW}$ or $f_\mathrm{CCW}$, depending on whether the observer is shifted into the moving FR according to a wave traveling clockwise (CW) or counterclockwise (CCW).

\section{Results}\label{sec:results}
\subsection{Comparison with standard dynamic mode decomposition}
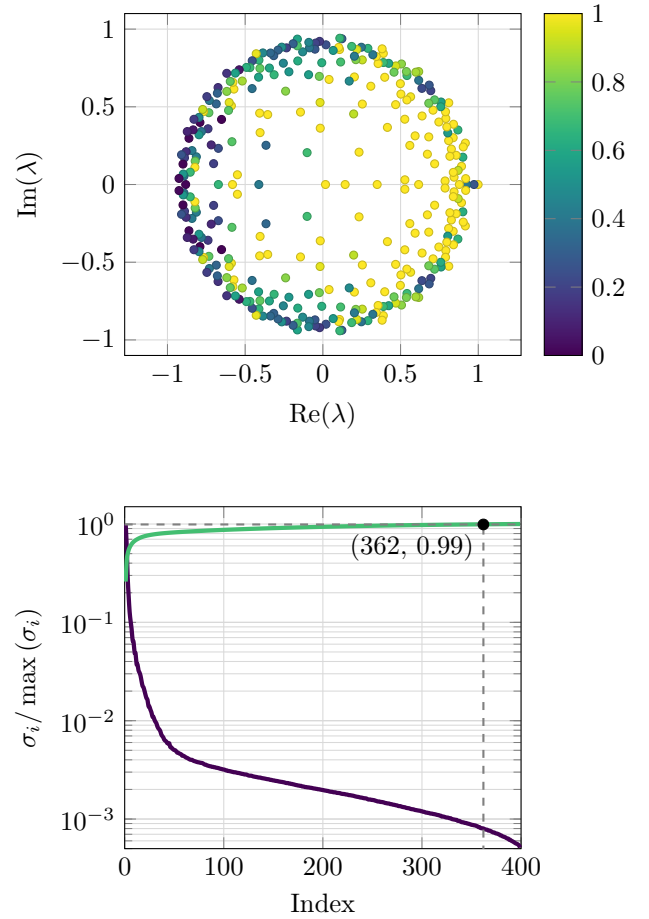
\begin{figure}[htpb]
	\centering
	\caption{The left panel shows the eigenvalue 
		distribution. The right panel shows the 
		normalized singular value decay where the 
		rank of truncation $k=362$ is indicated 
		to retain $99\%$ of the spatial energy 
	of the signal.}
	\pgfplotstableread[col sep = comma]{4_0_1_eigVal.csv} \eigenvalData 
\pgfplotstableread[col sep = comma]{4_0_2_sigmaDec.csv} \sigvalData 

\begin{tikzpicture}
    \begin{groupplot}[
        group style={group size=1 by 2, vertical sep=20mm}, 
        width=0.6\columnwidth, 
        xlabel = {$\operatorname{Re}(\lambda)$},
        ylabel = {$\operatorname{Im}(\lambda)$},
        scale only axis,
        axis on top,
        colormap/viridis,
        xlabel near ticks,
        ylabel near ticks,
	xmin = -1.1,
	xmax = 1.1,
	ymin = -1.1,
	ymax = 1.1,
        colorbar style={
            ylabel near ticks
        }
    ]
        \nextgroupplot[
            enlargelimits = false, 
            axis equal,
            grid=both,
            grid style={gray!30},
            colorbar,
            point meta min=0,
            point meta max=1,
            legend pos=north west,
            legend style={font=\scriptsize}
        ]
        \addplot[
            scatter,
            only marks, 
            mark=*, 
            mark size=1.5pt,
            scatter src=explicit,
        ] table[x index=0, y index=1, meta index=2] {\eigenvalData};
        
        \nextgroupplot[
            enlargelimits = false,
            ymode = log,
            xmin= 0,
            xmax= 400,
            ylabel = {$\sigma_i / \max \left(\sigma_i\right)$},
            xlabel = {Index},
            ymin= 5e-4,
            ymax= 1.5,
            grid=both,
            axis on top=false,
            grid style={gray!30},
            legend pos=north west,
            legend style={font=\scriptsize}
        ]
        \addplot[smooth, mark=none, ultra thick, color of colormap={0 of viridis}] table[x index=0, y index=1] {\sigvalData};
        \addplot[smooth, mark=none, ultra thick, color of colormap={700 of viridis}] table[x index=0, y index=2] {\sigvalData};
        \addplot[dashed, thick, gray] coordinates {(0, 0.99) (500, 0.99)};
        \addplot[dashed, thick, gray] coordinates {(362, 1e-10) (362, 1)};
        \addplot[only marks, mark=*, mark size=2pt] coordinates {(362, 0.99)} node[below left] {$(362,\, 0.99)$};
    \end{groupplot}
\end{tikzpicture}
	\label{fig:eigValSigDec}
\end{figure}

\begin{figure*}[htpb]
	\centering
	\caption{The left panel shows the original 
		data set and the middle panel the 
		reconstruction via the standard DMD method. 
		The right panel details the pointwise errors 
	between the data and its reconstruction.}
	\pgfplotstableread[col sep = comma]{4_0_3_errorVec.csv} \dataVecError 

\begin{tikzpicture}
    \begin{groupplot}[
        group style={group size=3 by 1, horizontal sep=15mm},
        width=0.2\textwidth,
        height = 5cm,
        scale only axis,
        axis on top,
        colormap/jet,
        xtick=\empty,
        ytick=\empty,
        xlabel near ticks,
        ylabel near ticks,
        colorbar style={
            yticklabel style={font=\scriptsize},
            ylabel near ticks
        }
        ]

        \nextgroupplot[
        enlargelimits = false, 
        ytick={0,16,32,48,64},
        yticklabels={600, 450, 300, 150, 0}, 
        xlabel = {$\theta$},
        ylabel = {$k\Delta t$},
        point meta min= 0 , point meta max= 100,
        colorbar style={
            yticklabel style={/pgf/number format/fixed, /pgf/number format/precision=2},
        },
        xtick={0,16,32,48,64},
        xticklabels={0, $\frac{\pi}{2}$, $\pi$, $3 \frac{\pi}{2}$, $2\pi$}, 
        title style={at={(0.1,-0.05)}, anchor=south, color = white},
        title = {(a)}
        ]
        \addplot graphics [xmin=0,xmax=64 ,ymin=0,ymax=64] {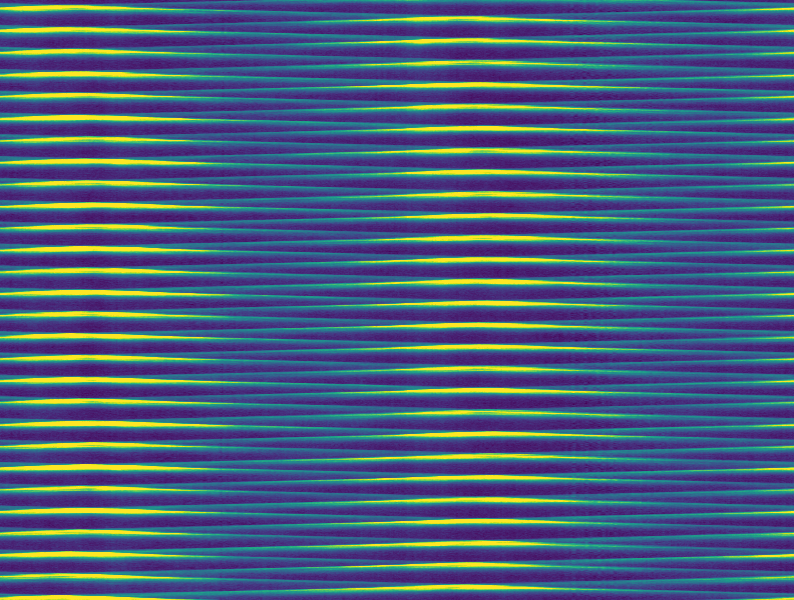};

        \nextgroupplot[
        enlargelimits = false, 
        xtick={0,16,32,48,64},
        xticklabels={0, $\frac{\pi}{2}$, $\pi$, $3 \frac{\pi}{2}$, $2\pi$}, 
        xlabel = {$\theta$},
        point meta min= 0  , point meta max=  100 ,
        colorbar style={
            yticklabel style={/pgf/number format/fixed, /pgf/number format/precision=2},
        },
                title style={at={(0.1,-0.05)}, anchor=south, color = white},
        title = {(b)}
        ]
        \addplot graphics [xmin=0,xmax=64 ,ymin=0,ymax=64] {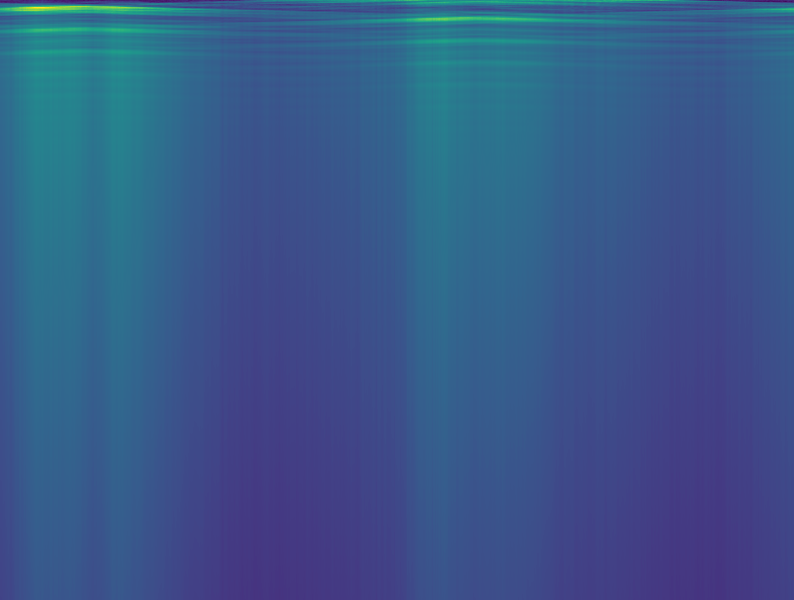};

        \nextgroupplot[
        enlargelimits = false, 
        xmode = log,
        xmin=0.1,
        xtick={0,0.1,1},
        y dir=reverse,
        xmax=1,
        xlabel = {$\frac{\left\|\boldsymbol{z}_k^{\mathrm{DMD}}-\boldsymbol{z}_k\right\|_2}{\left\|\boldsymbol{z}_k\right\|_2}$},
        width=0.1\textwidth,
        axis lines=left,
        tick align=outside,
        tick pos=left,
        grid=both,
        grid style={gray!30},
        legend pos=north west,
        legend style={font=\scriptsize},
                        title style={at={(0.2,-0.05)}, anchor=south, color = black},
        title = {(c)}
        ]
        \addplot[smooth, mark=none, teal, thick] table[x index=1, y index=0] {\dataVecError};
    \end{groupplot}
\end{tikzpicture}
	\label{fig:standardDMD}
\end{figure*}

In order to demonstrate the potential pitfalls RDE practitioners may encounter when analyzing counter-rotating detonation waves with the standard DMD (exactDMD) and to facilitate the necessity for the developed sequence,
the standard DMD is applied to the 2CR luminosity data in the following.
First,
the exactDMD procedure requires a rank of truncation,
which is selected according to the number of modes $k = 362$ required  to capture $99\%$ of the spatial energy.
The right panel in Figure~\ref{fig:eigValSigDec} shows the normalized decay of the singular values, 
exhibiting a slow decay typical for datasets that are corrupted with noise and transport-dominated dynamics such as traveling wave structures.
The left panel in Figure~\ref{fig:eigValSigDec} shows the distribution of the DMD eigenvalues in the complex plane,
where each marker is colored according to its weighting $|b_i|$ in the reconstruction.
First,
it becomes clear how the presence of measurement noise overestimates the real part of the DMD eigenvalues,
as nearly none of them lie on the unit circle,
in contrast to what is expected from the oscillatory dynamics of the counter-rotating detonation waves.
Second,
it shows that an energy-based heuristic for selecting the rank of truncation also overestimates $k$,
which is indicated by the presence of erroneous DMD eigenvalues.
Manually discriminating the DMD eigenvalues into those that carry actual dynamical information from those corrupted by noise is cumbersome and prone to error.
Furthermore,
this discrimination is even harder for the 2CR case since the primary and secondary wave DMD eigenvalues are indistinguishable in the spectrum.
At the same time,
it demonstrates the benefits of changing the frame of reference,
as the DMD eigenvalues and their harmonics degenerate to near-stationary eigenvalues.
Finally,
the standing wave pattern and the periodically occurring nonlinear interaction of the counter-rotating detonation waves justify 
the use of EDMD in combination with time delay embedding.
The interplay of the discussed shortcomings of the standard DMD can be seen in Figure~\ref{fig:standardDMD}(c),
which shows the relative error between the data and its reconstruction shown in Figure~\ref{fig:standardDMD}(a) and (b) at a given time $t$.

\subsection{Algorithmic sequence for the analysis of counter-rotating detonation waves}
\begin{figure*}[htpb]
	\centering
	\def \svgwidth{0.9\textwidth}
	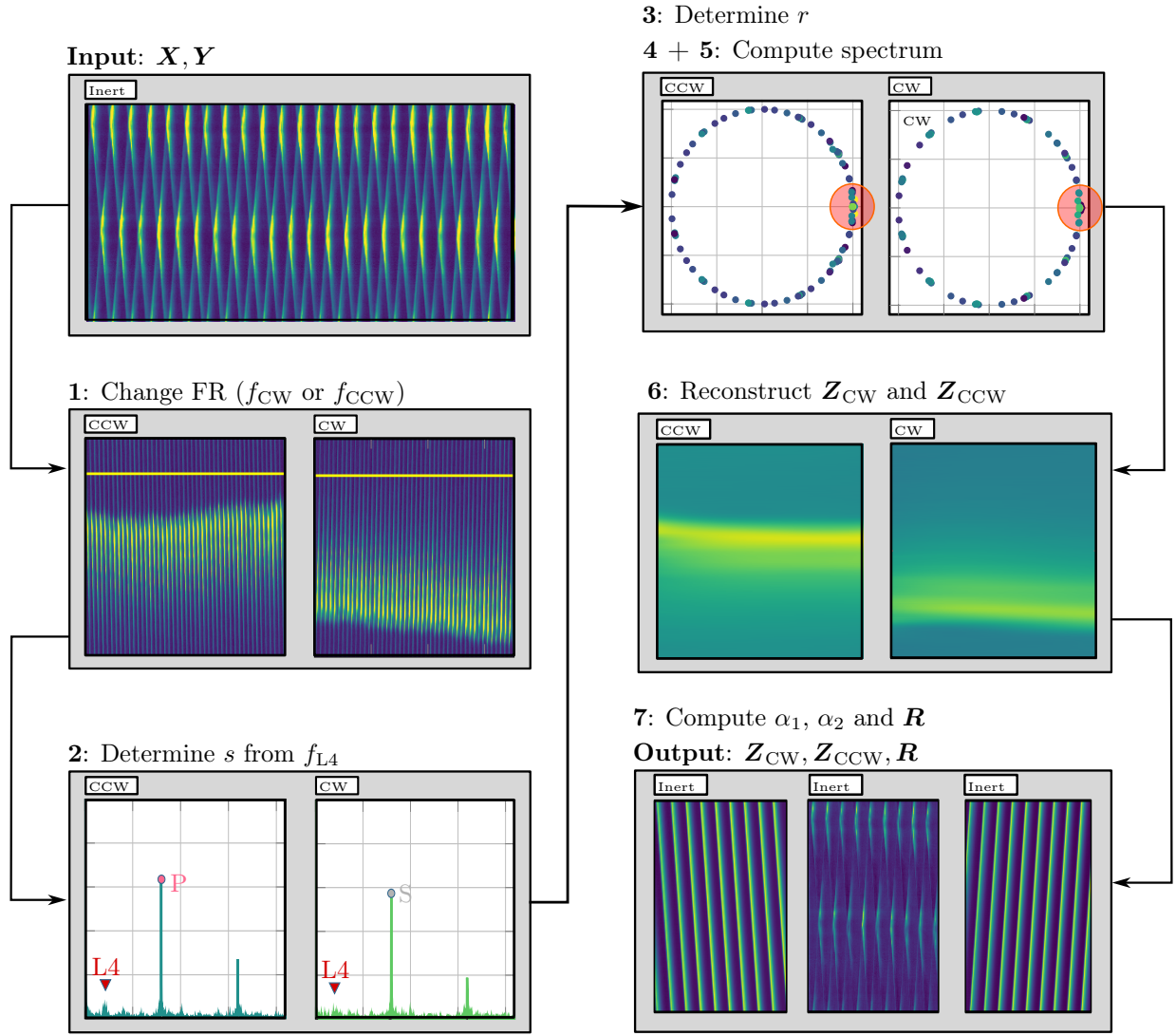
	\caption{The flow chart illustrates the individual steps of the algorithmic description of the sequence. 
		In each step the FR is highlighted.
		Further,
		in the panels that correspond to steps 3, 4 and 5 red circles 
	indicate the eigenvalues used for reconstructing the travelling wave in the corresponding FR respectively.}
	\label{fig:scheme}
\end{figure*}

We propose  an  \algoSeq that combines the methods discussed in section~\ref{sec:methods} into a series  of DMD-like algorithms and data-processing steps,
which are summarized in Algorithm~\ref{alg:algo-seq}  as pseudo code and are illustrated in Figure~\ref{fig:scheme} as a flow chart. 
Among established methods such as ExactDMD and OptDMD, a central element of the sequence is the physics informed estimation of the number of time delay embeddings based on the periodicity of the dynamics.
This allows for accurate DMD models to be built  in the according moving FR,
which in turn enables a precise decomposition of the operating mode under consideration into datasets that capture the traveling waves and a nonlinear remainder.    
First,
in Step~\ref{step:change-fr},
the FR in which the data is represented is changed from the inertial to the moving FR according to the frequency $f_\mathrm{CW}$ or $f_\mathrm{CCW}$ of the wave traveling clockwise (CW) or counterclockwise (CCW),
respectively.  
In Step~\ref{step:determine-s},
the number of time-delay embeddings of the full state, $s$, is determined from the lowest dynamically relevant frequency, $f_\mathrm{Base}$.  
Subsequently,
in Step~\ref{step:determine-r},
the rank of truncation, $r$,
is computed via the optimal-threshold criterion for the extended data matrices.  
In Step~\ref{step:exact-dmd},
the DMD eigenvalues and DMD modes are computed via the ExactDMD algorithm.  
Using the results of ExactDMD as an initial guess,
in Step~\ref{step:opt-dmd} OptDMD is applied to remove the bias introduced by noise.
In Step~\ref{step:reconstruct} the datasets describing the traveling wave structures are then reconstructed via
\begin{equation}\label{eq:dmd-reconstruct-feature}
	\boldsymbol{z}_k^\mathrm{CCW/CW}  \approx \sum_{j=1}^t\left(\boldsymbol{B}^{\top} \boldsymbol{\phi}_j\right) \lambda_j^k b_j
\end{equation}
by selecting the according subset of the computed spectrum.     
Finally,
in Step~\ref{step:lsq} the scaling coefficients for the reconstructed subset are computed by solving the least square problem 
\begin{equation} \label{eq:lsq-coeffalpha}
	\begin{bmatrix} 
		\alpha_1 \\ 
		\alpha_2 
	\end{bmatrix}
	= \arg \min _{ 
		\tilde{\alpha}_1, \tilde{\alpha}_2
	}\left\|
	\boldsymbol{Z}_{\text{train}}
	- \left(\tilde{\alpha}_1 \boldsymbol{Z}_{\mathrm{CW}}+\tilde{\alpha}_2 \boldsymbol{Z}_{\mathrm{CCW}}\right)
	\right\|^2_F.
\end{equation} 
Subsequently,
the remainder containing the nonlinear interactions is computed as
\begin{equation}\label{eq:remainder}
	\boldsymbol{R}=\boldsymbol{Z}-\alpha_1 \boldsymbol{Z}_{\mathrm{CW}}-\alpha_2 \boldsymbol{Z}_{\mathrm{CCW}}.
\end{equation}
\begin{algorithm}
	\caption{Algorithmic sequence for the analysis of counter-rotating detonation waves}
	\label{alg:algo-seq}
	\KwIn{Snapshot matrices $\boldsymbol{X}, \boldsymbol{Y} \in \mathbb{R}^{ns \times m}$}
	Change FR with respect to $f_\mathrm{CW}$ or $f_\mathrm{CCW}$\; \label{step:change-fr}
	Determine number of time-delay embeddings $s$\; \label{step:determine-s}
	Determine rank of truncation $r$\; \label{step:determine-r}
	Compute spectrum of $\tilde{\boldsymbol{K}}_\mathrm{DMD}$ via ExactDMD\; \label{step:exact-dmd}
	Refine spectrum of $\tilde{\boldsymbol{K}}_\mathrm{DMD}$ via OptDMD\; \label{step:opt-dmd} 
	Reconstruct $\boldsymbol{Z}_{\mathrm{CW}}$ and $\boldsymbol{Z}_{\mathrm{CCW}}$\; \label{step:reconstruct}
	Compute for scaling $\alpha_1$ and $\alpha_2$ and compute remainder $\boldsymbol{R}$.\; \label{step:lsq}
	\KwOut{$\boldsymbol{Z}_{\mathrm{CW}}$, $\boldsymbol{Z}_{\mathrm{CCW}}$ and $\boldsymbol{R}$} 
\end{algorithm}

\subsection{Application of the sequence}
Across the operating modes,
the number of time steps between collisions of counter-rotating waves ranges from 5 to 10 time steps.  
Therefore,
from each dataset containing up to 5000 snapshots,
a subset of 500 snapshots is chosen,
as it captures a sufficient number of collisions of counter-rotating waves while keeping the computational cost manageable.  
The subsets are chosen in such a way that the traveling waves appear approximately stationary in the CW and CCW moving FR.

\begin{figure}[htpb]
	\centering
	\def \svgwidth{\columnwidth}
	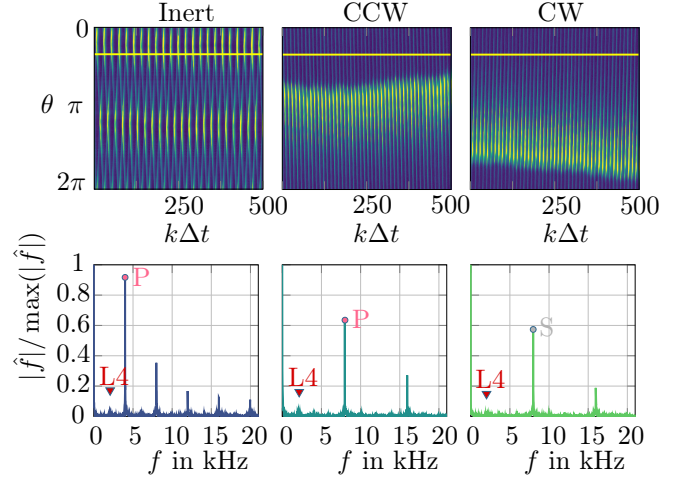
	\caption{ Overview of the selected subsets for the 2CR operating mode, spanning $m_{\mathrm{train}} = 500$ time steps and the FFT taken at the azimuthal position $\theta_\mathrm{PCB} = 60\SI{}{\degree}$ in the inertial FR (left column), the CCW moving FR (middle column), and the CW moving FR (right column).} 
	\label{fig:4_1_2CR_FR}
\end{figure}

The top row in Figure~\ref{fig:4_1_2CR_FR} shows an $\theta$- $t$ diagram of  dataset representing the 2CR operating mode in the inertial FR in the left column,
in the moving FR associated with the wave traveling CCW in the middle column,
and in the moving FR associated with the wave traveling CW in the right column.  
The frequencies $f_\mathrm{CW}$ and $f_\mathrm{CCW}$ which are used as reference for changing the FR are obtained from an FFT of the signal in the inertial FR and are determined as $f_\mathrm{CW} = f_\mathrm{CCW}=\SI{3989}{\hertz}$. 

The bottom row of Figure~\ref{fig:4_1_2CR_FR} shows the FFT of the signals obtained from the luminosity probes at the azimuthal positions $\theta_{\mathrm{PCB}}$ indicated by the yellow lines in the top row.
In the left column the spectrum is shown in the inertial FR as in Section~\ref{sec:analysis-operating-mode}.
The middle column and the right column show the spectrum sampled at the same position but in the CCW moving FR and the CW moving  FR respectively.
It shows how the frequencies of the counter-rotating waves appear accelerated in the moving FR.
The peaks that belong to the wave used as reference for the transformation vanish in the moving FR.
The counter-rotating wave then appears faster in the moving FR by the reference frequency.
For example in the 2CR case the peaks of the primary wave in the CCW moving FR satisfy
\begin{math}
	f_{\mathrm{P,CCW}} = f_{\mathrm{CCW}} + f_{\mathrm{P}}
\end{math}
and the peaks of the secondary wave vanish.
Conversely in the CW moving FR the peaks of the secondary wave satisfy
\begin{math}
	f_{\mathrm{S,CW}} = f_{\mathrm{CW}} + f_{\mathrm{S}} =  f_{\mathrm{P,CCW}}
\end{math}
and the frequencies associated with the peaks of the primary wave vanish.
Furthermore, it  is observed that the frequency of the longitudinal mode L4 is preserved under the transformation.
With the spectral information at hand it is now possible to estimate the number of time delay embeddings $s$ based on the periodicity of the data.  
To this end, the lowest dynamically relevant frequency $f_{\mathrm{Base}}$ is determined from the FFT for each moving FR. The number of time delay embeddings of the full state is then computed as
\begin{equation}
	\label{eq:stacks}
	s = \left\lfloor \frac{T_{\mathrm{Base}}}{\Delta t} \right\rfloor
	= \left\lfloor \frac{1}{f_{\mathrm{Base}}\,\Delta t} \right\rfloor.
\end{equation}

Across all considered operating modes in this dataset, the lowest dynamically relevant frequency $f_{\mathrm{Base}}$ is identified as the frequency corresponding to the peak of the L4 mode.
Since the frequency of the longitudinal mode L4 remains unchanged under the change of FR it is identified as the lowest relevant frequency $f_{\mathrm{Base}}$ for both the CW and CCW moving FR.

Given the number of time delay embeddings the required number of modes $r$ can be determined as discussed in section \ref{sec:optimal-treshold}.  

\begin{figure*}[htpb]
	\centering
	\def \svgwidth{\textwidth}
	\scriptsize{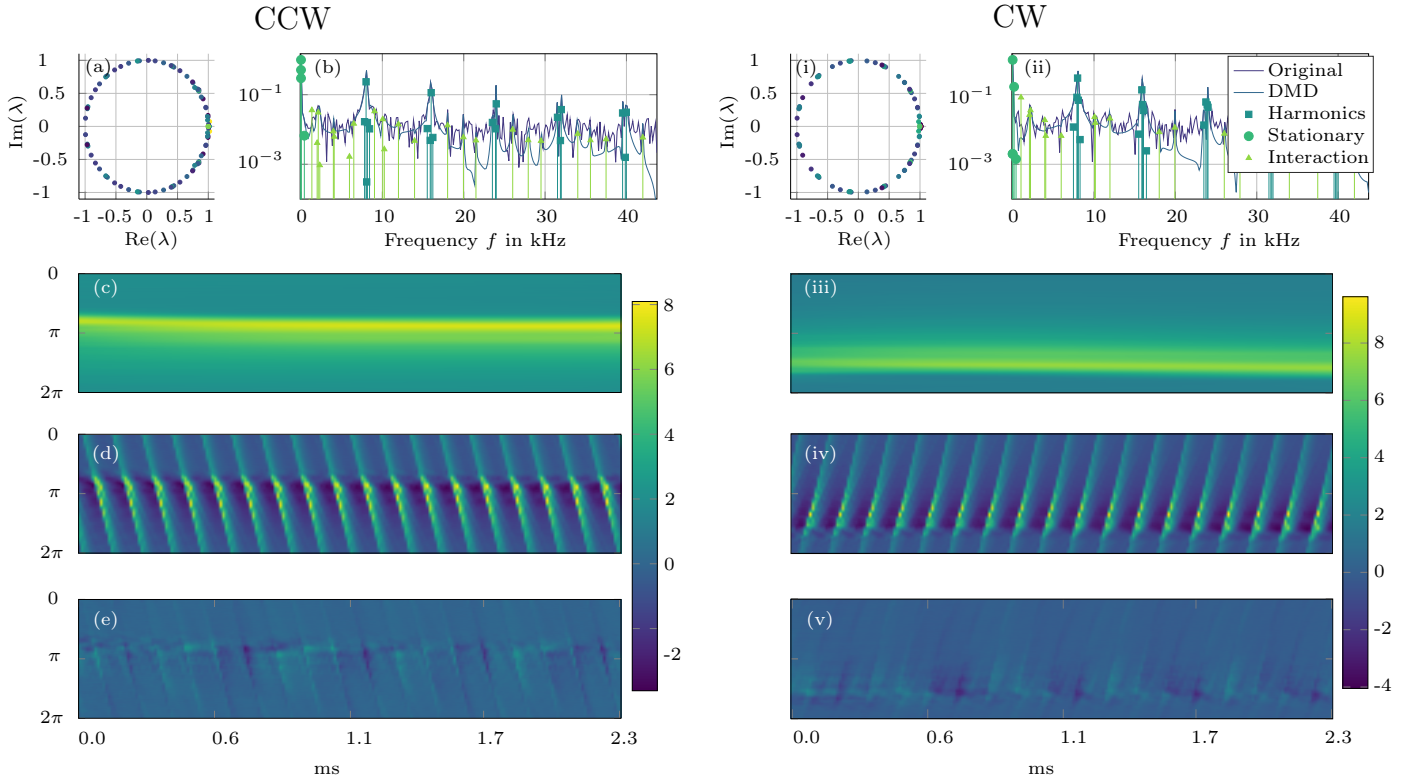}
	\caption{(a) and (i) show the discrete spectrum of the operator $K_\mathrm{DMD}$ built on the data  in the CCW and CW moving FR,  panel (b) and (ii)  compare the FFT of the signal and its reconstruction at $\theta_\mathrm{PCB} = 60\SI{}{\degree}$ in the CCW and CW moving FR. Markers indicate the assignment of the frequencies to the dynamical features. The bottom panels (c), (d), (e) and (iii), (iv), (v) show their reconstruction via DMD in the CCW and CW moving FR.}
	\label{fig:4_3_1_DMD}
\end{figure*}

Finally, for the 2CR operating mode, the required parameters $r$ and $s$ for the sequence are determined as $s_\mathrm{Base} = 43$ and  the number of modes required are $r_{\mathrm{CCW},s}=90$ and $r_{\mathrm{CW},s}=86$.
The sensitivity of the DMD models with respect to the parameters $s$ and  $r$ in terms of the reconstruction error and spectral agreement is discussed in Appendix~\ref{sec:sens-ana-appendix}.
With the parameters at hand we proceeded to build the DMD models in the CCW and CW moving FR.
Figure~\ref{fig:4_3_1_DMD} shows the  results of the modeling in the CCW  and CW moving FR.
Figure~\ref{fig:4_3_1_DMD} (a) and (i) show the DMD eigenvalues after post-processing via the OptDMD procedure in the CCW moving and CW moving FR, respectively.
It shows that all identified DMD eigenvalues lie approximately on the unit circle,
indicating that most of the bias introduced by noise has been successfully removed.

We then proceed to relate the resulting DMD eigenvalues to the dynamical features of the system in the selected FR through the relation $\omega_j = \frac{\ln\left( \lambda_j \right)}{ \Delta t }$ where $\omega_j$ denotes the  time  continuous eigenvalue. Consequently, the frequency associated with the time continuous eigenvalue is given by
\begin{equation}\label{eq:freq-time}
	f_j = \frac{1}{2\pi}\mathrm{Im}(\omega_j).
\end{equation}
Therefore,
by comparing the frequencies associated with the DMD eigenvalues computed as in~\eqref{eq:freq-time} with the frequencies obtained from the FFT of the luminosity probe,  we can identify which DMD modes correspond to which dynamical features.
Figure~\ref{fig:4_3_1_DMD} (b) and (ii)  compares the FFT of the original signal with the FFT of the reconstructed signal sampled at $\theta_\mathrm{PCB}$ in the moving FR respectively.
In general, it can be seen that across all identified modes, there is strong agreement between the frequencies associated with the DMD eigenvalues and the peaks of the FFT spectrum. 
This agreement allows for a direct assignment of DMD eigenvalues and DMD modes to physical dynamical features.
Figure~\ref{fig:4_3_1_DMD} (c), (d), (e) and (iii), (iv), (v) show the reconstruction of the individual  dynamical features such as the primary wave (P), secondary wave (S), interactive modes as in~\eqref{eq:dmd-reconstruct-feature} based on the assignment of DMD eigenvalues to dynamical features.

In summary, by changing the FR twice, it becomes possible to  approximately separate the traveling wave structures (P) and (S) from the remaining dynamics of the operating mode.  
Therefore, the traveling waves are first transformed back into the inertial FR. 
Subsequently, the scaling coefficients  $\alpha_1$ and  $\alpha_2$ are computed as in \eqref{eq:lsq-coeffalpha} and the remainder $\boldsymbol{R}$ is computed as in \eqref{eq:remainder}. 

\begin{figure*}[htbp]
	\centering
	\def \svgwidth{.8\textwidth}
	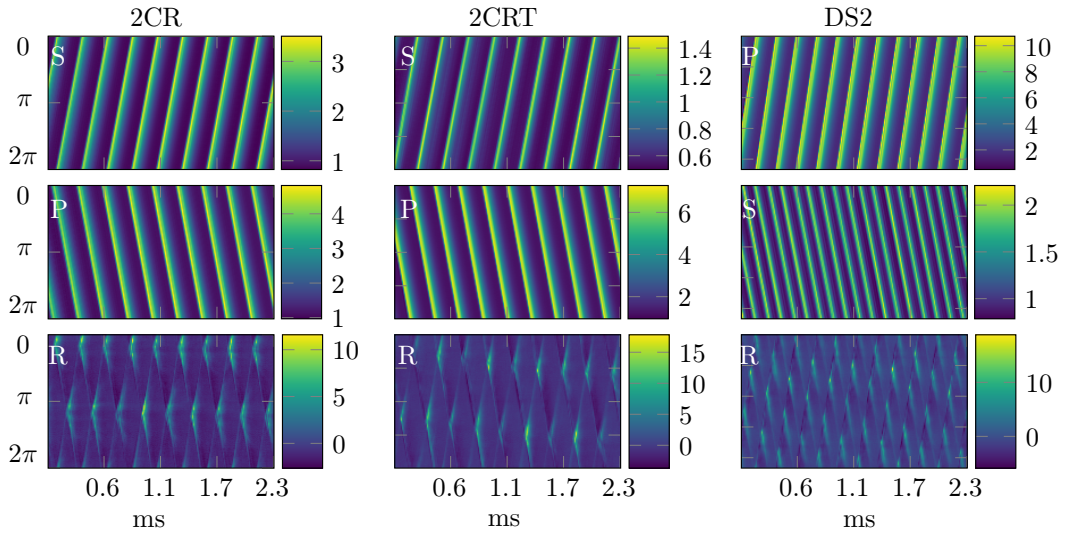
	\caption{The figure shows the temporal evolution of the identified traveling wave structures for the different operating modes. The top  row shows the waves traveling in CCW direction, the middle row shows the waves traveling CW and the bottom row shows the nonlinear remainder for each operating modes respectively.}
	\label{fig:4_4_extractedFeatures}
\end{figure*}

An overview of the output generated by the algorithmic sequence  \ref{alg:algo-seq} for other operating modes is shown in Figure~\ref{fig:4_4_extractedFeatures}. It shows the discrimination of the operating mode under consideration into the traveling wave structures in top and middle row (labeled P and S). The bottom row in Figure~\ref{fig:4_4_extractedFeatures}
then shows the nonlinear remainder $\boldsymbol{R}$.
\subsection{Nonlinear wave amplification}\label{sec:comp-interaction}
\begin{figure*}[ht]
	\centering
	\def \svgwidth{\textwidth}
	\footnotesize{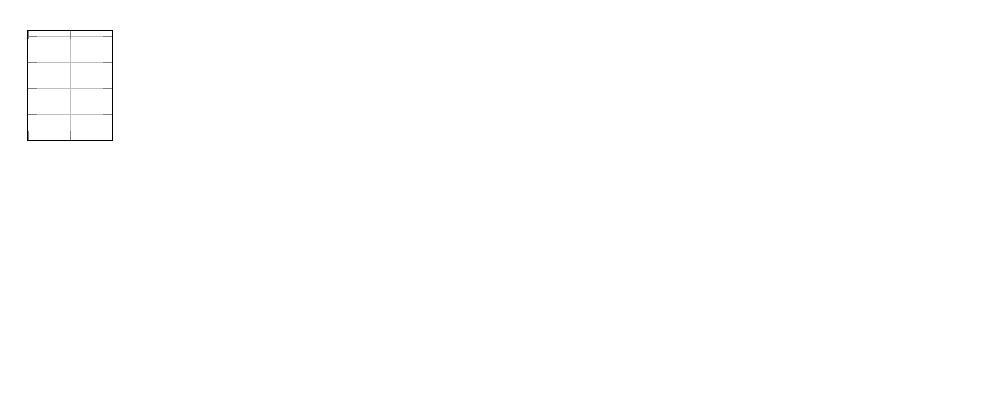}
	\caption{The top row shows the full luminosity signal sampled at times around the collision.  
		The middle row shows the CCW and CW traveling waves and their crests indicated by red and yellow triangular markers respectively.  
	The bottom row shows the remainder signal $\boldsymbol{R}$ after subtracting the scaled traveling waves, where blue triangular markers indicate maximum luminosity.}
	\label{fig:4_5_nonLinAmplification}
\end{figure*}

With a precise estimate of the traveling waves at hand,
the next step is to quantify the nonlinear interaction present in the operating modes. 
Figure~\ref{fig:4_5_nonLinAmplification} shows a sequence of snapshots in temporal proximity to the collision of the two counter-rotating waves for each of the operating modes.  
The top row shows the snapshots sampled from the original data and the middle row shows snapshots from the reconstructions of the primary and secondary waves only.
Furthermore,
the crest of the wave traveling CCW is marked with a red triangular marker and the crest of the wave traveling CW is marked with a yellow triangular marker.  
The bottom row displays snapshots taken from the remainder $\boldsymbol{R}$.
The blue markers indicate the maximum of each snapshot.
Furthermore,
it is indicated how much the nonlinear signal contributes to the total luminosity.  
It illustrates how the nonlinear amplification becomes strongest at the collision of the two counter-rotating waves and then decays as the crests of the traveling waves move apart from each other.   

For the 2CR operating mode,
at the point of wave collision, which is indicated by the red line,
\qty{54.47}{\percent} of the luminosity is made up by the nonlinear amplification.
For the 2CRT operating mode,
at the point of collision \qty{61.57}{\percent} of the luminosity is made up by nonlinear amplification.
Furthermore,
for the DS2 operating mode,
at the point of collision at one of the two waves crest in the,
indicated by the red line,
the luminosity of the signal is  made up by \qty{55.79}{\percent} of nonlinear amplification.

\subsection{Comparison of the dynamics of the nonlinear interactions}\label{sec:comp-dyn}
\begin{figure*}[htpb]
	\centering
	\pgfplotstableread[col sep = comma]{4_7_trunkModeSens_specAgree_CCW_2CR.csv}  \tmSACCWIICR
\pgfplotstableread[col sep = comma]{4_7_trunkModeSens_specAgree_CCW_2CRT.csv} \tmSACCWIICRT
\pgfplotstableread[col sep = comma]{4_7_trunkModeSens_specAgree_CCW_DS2.csv}  \tmSACCWDS
\pgfplotstableread[col sep = comma]{4_7_trunkModeSens_specAgree_CW_2CR.csv}   \tmSACWIICR
\pgfplotstableread[col sep = comma]{4_7_trunkModeSens_specAgree_CW_2CRT.csv}  \tmSACWIICRT
\pgfplotstableread[col sep = comma]{4_5_0_nonlinearRemainder_2CR.csv} \dataRemainderTwoCR
\pgfplotstableread[col sep = comma]{4_5_1_nonlinearRemainder_2CRT.csv} \dataRemainderTwoCRT
\pgfplotstableread[col sep = comma]{4_5_2_nonlinearRemainder_DS2.csv} \dataRemainderDS

\pgfplotscreateplotcyclelist{three_col_dense}{
    {color of colormap={0 of viridis}, mark=*, line width=1.0pt,
        solid,
        mark repeat=2, mark phase=0, mark options={draw=black}},
    {color of colormap={300 of viridis}, mark=square*, line width=1.0pt,
        dashed,
        mark repeat=2, mark phase=1, mark options={draw=black}},
    {color of colormap={600 of viridis}, mark=triangle*, line width=1.0pt,
        dotted,
        mark repeat=2, mark phase=2, mark options={draw=black}},
}

\begin{tikzpicture}
    \begin{groupplot}[
        group style={group size=1 by 3, vertical sep=12mm},
        height=2.2cm,
        width = 0.8\textwidth, 
        scale only axis,
        enlargelimits=false,
        tick pos=left,
        grid=both,
        grid style={gray!30},
        xlabel near ticks,
        ylabel near ticks,
        clip=false,
        legend style={
            at={(1.01, 0.5)}, 
            anchor=west,      
            cells={anchor=west}
        }
    ]

    \nextgroupplot[
        ylabel={$L$},
        title={2CR},
        cycle list name = three_col_dense,
    ]
    \addplot table[x index=0, y index=1] {\dataRemainderTwoCR};
    \addplot table[x index=0, y index=2] {\dataRemainderTwoCR};
    \legend{$\boldsymbol{R}_P$, $\boldsymbol{R}_S$}
    \node[anchor=south west, outer sep=2pt] at (rel axis cs:0.0, 1.0) {\textbf{(a)}};

    \nextgroupplot[
        ylabel={$L$},
        title={2CRT},
        cycle list name = three_col_dense,
    ]
    \addplot table[x index=0, y index=1] {\dataRemainderTwoCRT};
    \addplot table[x index=0, y index=2] {\dataRemainderTwoCRT};
    \legend{$\boldsymbol{R}_P$, $\boldsymbol{R}_S$}
    \node[anchor=south west, outer sep=2pt] at (rel axis cs:0.0, 1.0) {\textbf{(b)}};

    \nextgroupplot[
        ylabel={$L$},
        xlabel={$k\Delta t$},
        title={DS2},
        cycle list name = three_col_dense,
    ]
    \addplot table[x index=0, y index=1] {\dataRemainderDS};
    \addplot table[x index=0, y index=2] {\dataRemainderDS};
    \addplot table[x index=0, y index=3] {\dataRemainderDS};
    \legend{$\boldsymbol{R}_P$, $\boldsymbol{R}_{S,i}$, $\boldsymbol{R}_{S,ii}$}
    \node[anchor=south west, outer sep=2pt] at (rel axis cs:0.0, 1.0) {\textbf{(c)}};
\end{groupplot}
\end{tikzpicture}
	\caption{The figure compares the amplitude of the nonlinear remainder  sampled at the azimuthal position of the crest of the counter-rotating waves.}
	\label{fig:4_6_qualitativeInteraction}
\end{figure*}
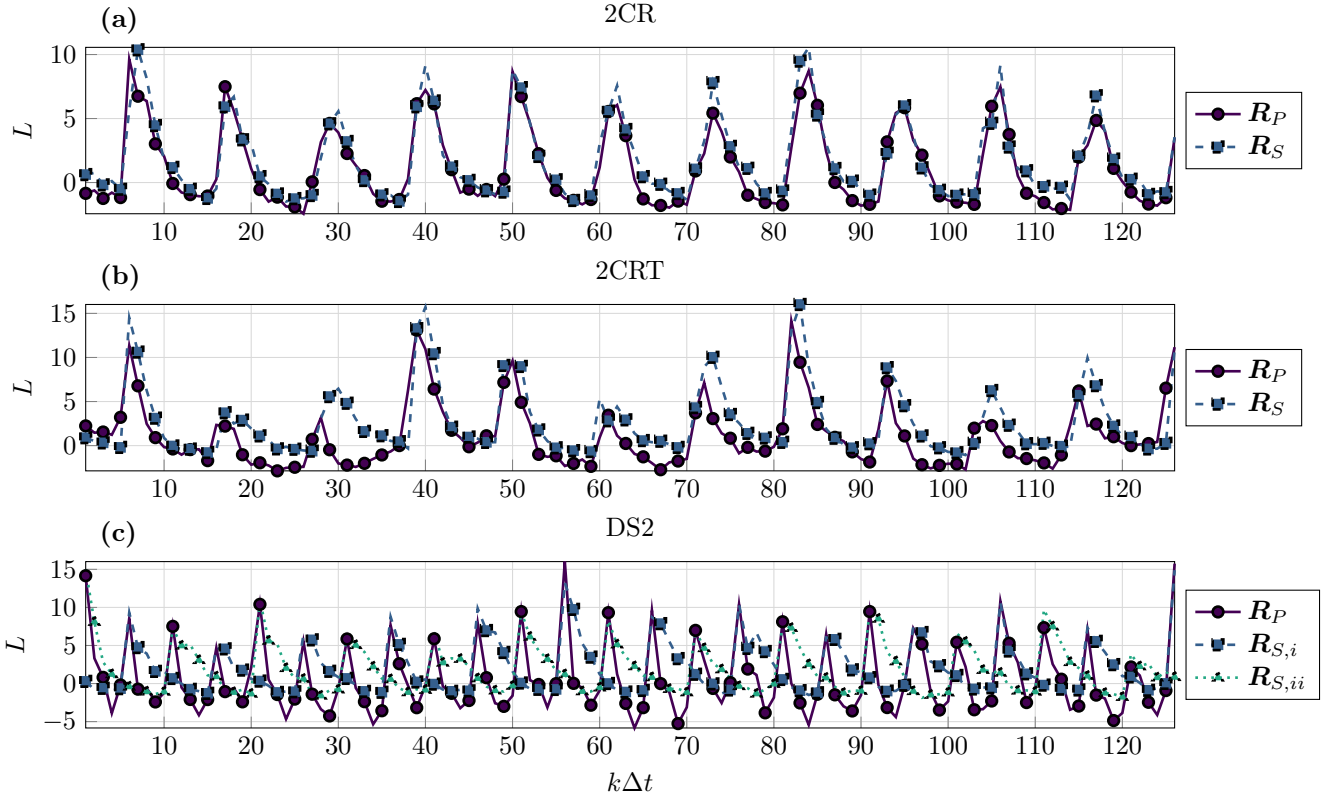

The analysis  outlined in Section~\ref{sec:comp-interaction} suggests that the strength of the nonlinear amplification at the collision point of counter-rotating waves depends on the distance of the crests from the point of collision and the operating mode under consideration.
In order to further characterize the nonlinear interaction, the remainder $\boldsymbol{R}$ is sampled at the azimuthal position $\theta_\mathrm{max}$,
given by  locations of the crests of the counter-rotating waves under consideration (as indicated by the yellow and red triangular markers in Figure~\ref{fig:4_5_nonLinAmplification}).
By doing so,
the change in amplitude of the nonlinear remainder at the crests position is traced in dependency of the distance to the collision.
Figure~\ref{fig:4_6_qualitativeInteraction}\,(a) shows that for the 2CR case,
the remainder of the crest of the primary ($\boldsymbol{R}_\mathrm{P}$) and of the secondary ($\boldsymbol{R}_\mathrm{S}$) waves behave similarly.
At the points of collision,
the luminosity increases sharply and then decays gradually until the next collision occurs.
Figure~\ref{fig:4_6_qualitativeInteraction}\,(b) shows the remainder luminosity at the crest of the primary ($\boldsymbol{R}_\mathrm{P}$) and secondary ($\boldsymbol{R}_\mathrm{S}$) wave for the 2CRT case.
Similar to the 2CR case,
the luminosity increases sharply at the point of collision.
However,
after the collision, 
the remainder associated with the crest of the primary wave decays faster than that of the secondary wave.
This suggests that the interaction affects the primary and secondary waves differently.
Finally,
Figure~\ref{fig:4_6_qualitativeInteraction}\,(c) shows the remainder luminosity at the wave crest of the primary wave ($\boldsymbol{R}_\mathrm{P}$) and the two crests of the set of counter-rotating waves ($\boldsymbol{R}_\mathrm{S,i}$ and $\boldsymbol{R}_\mathrm{S,ii}$). 
Once more, at every collision between the primary wave and one of the counter-rotating waves,
the remainder increases sharply.
Furthermore,
as seen in the 2CRT case,
the remainder associated with the primary wave decays significantly faster than that of the counter-rotating waves.
The fact that the counter-rotating waves decay more slowly in both the 2CRT and DS2 cases supports the conclusion that the primary wave plays a stabilizing role in the interaction with the counter-rotating waves.

Ultimately, we can conclude that the nature of the interactions between the counter-rotating waves in each of these modes behaves differently. Unfortunately, the natural luminosity of the flame measured in these datasets is limited in what we can interpret from the data, except to qualitatively say that the extent and amount of luminosity will more or less correlate with the temperature of the products. From this, we can interpret that the 2CR waves (primary and secondary) behave more or less indistinguishably from each other. However, slightly slower decay rate of the luminosity in secondary wave in the 2CRT case and the significantly slower decay rate of the luminosity in the secondary waves in the DS2 case indicates that the collision of the secondary wave with the primary wave provides a periodic support to the wave. The nonlinear residual attached to each of the secondary waves indicates that it is being re-energized by each collision. This is especially prevalent in the DS2 case where the primary wave rapidly decays back to its unperturbed state after each collision, where the secondary waves preserve this nonlinearity for a much longer time. With this dataset, it is hard to estimate the net effect on the overall performance of the combustor. However, the increased presence of the nonlinearity in the collisions of the primary and secondary waves in the DS2 case may be associated to a stronger coupling of the detonation wave. That said, further experiments would be required to provide greater insights into the actual heat release process in the waves, e.g. through measurements of the OH* chemiluminescence as a proxy for heat release.

\section{Conclusion and outlook}\label{sec:discussion}
In this work,
we developed an \algoSeq to  characterize and quantify the dynamics of  traveling detonation waves in a RDC. 

First,
a data post-processing routine was introduced to reduce the dimensionality and to account for experimental influences in the luminosity data generated from high-speed video imaging of the aft end of the RDC.  
Subsequently,
Koopman operator theory was presented as a framework to describe the nonlinear dynamics.  
It was shown that the EDMD,
combined with a dictionary constructed from time-delay embeddings of the full state,
provides meaningful finite-dimensional approximations of the Koopman operator and is suited for the RDC data.  
In order to compute those approximations,
the ExactDMD was introduced as an algorithm that can take advantage of low-dimensional structures in the data.  
Furthermore,
to address the influence of sensor noise in the measurements,
the OptDMD method was introduced as a post-processing step for the results generated by ExactDMD.  
Finally,
an optimal threshold criterion and an observer change based on the traveling wave speed were introduced.  
The individual steps were then assembled into an algorithmic sequence for analyzing detonation waves.

Subsequently,
it was shown that the developed algorithmic sequence is capable of extracting the traveling wave structure from the  data describing the different operating modes.   
Moreover,
by subtracting the traveling waves from the full signal,
a remainder containing the nonlinear amplification was obtained.  

It was possible to quantify the contributions to the luminosity signal due to nonlinear amplification and due to linear interference.  
It was shown that the contribution due to nonlinear amplification depends on the operating mode,
which aligns with the results in \cite{barnouinRDC25}.

In addition,
the nonlinear interaction of counter-rotating waves could be further characterized by sampling the remainder at the azimuthal position of the wave crests.  
In particular,
it was shown that for the 2CRT and DS2 operating modes,
the primary and secondary wave structures are affected differently.  
Whereas the primary wave decays faster,
the secondary decays much slower,
which indicates that the primary wave stabilizes the secondary wave structures.

The precise separation of the traveling wave structures from the detonation wave opens several directions for future work.  
As an intermediate next step,
the algorithmic sequence should be extended such that the remainder itself is subject to DMD analysis.  
This may help to separate the longitudinal modes that are still present in the data and could offer additional insight into the structure of the interaction terms.

Another promising extension would be to use the algorithmic sequence on the pressure data.  
Although the spatial resolution is lower than for the luminosity data, the temporal resolution is higher and may still allow for meaningful DMD models to be constructed.

Another popular approach for system identification is the Sparse Identification of Nonlinear Dynamical Systems (SINDy) method \cite{bruntonSindy16},
which retrieves governing equations from a symbolic regression performed on measured data.  
The work by \cite{forootaniSindy24} extended this idea to the identification of partial differential equations using a greedy method,
which could be applied to the RDC data as well.

A further point concerns the choice of dictionary $\mathcal{D}$ required for the EDMD method.  
Throughout this work,
a dictionary based on time-delay embedding was exclusively considered,
which proved to be useful.  
However,
approaches from machine learning,
such as deep learning could be leveraged to find alternative embeddings in the context of Koopman operator theory as done in the work by \cite{goyalDeepLearn23, luschNNDMD18}.
The data that support the findings of this study are openly available  at \cite{oexleHSV2026}.

\appendix
\section{Koopman operator theory} \label{sec:koop-operator-theory-appendix}
A nonlinear dynamical system may  be written in its general form as  
\begin{equation}\label{eq:dyntimecont}
	\frac{\mathrm{d}}{\mathrm{d} t} {\boldsymbol{z}}(t)=\boldsymbol{f}({\boldsymbol{z}}(t))
\end{equation}
where $\boldsymbol{f}$ is a  (possibly unknown) vector field that describes the evolution of  the state  $\boldsymbol{z}\left(t \right) \in \mathcal{M}$  on the state space $\mathcal{M} \subseteq \mathbb{R}^n$.

The nonlinear flow map $\boldsymbol{F}\colon \mathbb{R}^n \rightarrow \mathbb{R}^n$ is associated with the time-continuous dynamical system \eqref{eq:dyntimecont} through 
\begin{equation}\label{eq:dyntimedisc}
	\boldsymbol{F}(\boldsymbol{z}_k)= \boldsymbol{z}_k + \int_{t_k}^{t_{k+1}} \boldsymbol{f}({\boldsymbol{z}}(\tau)) \mathrm{d}\tau = \boldsymbol{z}_{k+1}.
\end{equation}

Instead of acting on the state $\boldsymbol{z}$,
the Koopman operator acts linearly \cite{bruntonBook21}  on a scalar observable function $g$ as
\begin{equation}
	\mathcal{K} g\left(\boldsymbol{z}_k\right)=g\left(\boldsymbol{F}\left(\boldsymbol{z}_k\right)\right)=g\left({\boldsymbol{z}}_{k+1}\right),
\end{equation} where  $g  \in \mathcal{F}$ with $g\colon\mathcal{M}\to \mathbb{R}$ and $\mathcal{F}$ denotes an infinite-dimensional function space (e.g. $L^2(\mathcal{M})$).

\begin{figure}[htpb]
	\centering
	\def \svgwidth{\columnwidth}
	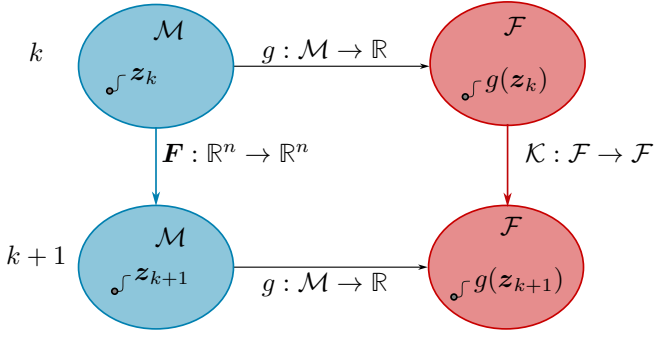
	\caption{The figure shows the action of the Koopman operator on the scalar valued observable function.}
	\label{fig:koopman_action}
\end{figure}

Figure~\ref{fig:koopman_action} compares these two possibilities of describing the evolution of a dynamical system. 
Either  the finite-dimensional state is  advanced  via the nonlinear flow map $\boldsymbol{F}\colon \mathbb{R}^n \rightarrow \mathbb{R}^n$ or  the  function of the state,
i.e.,
the observable function,
is advanced via the infinite-dimensional Koopman operator $\mathcal{K}\colon \mathcal{F} \rightarrow \mathcal{F}$.

A set of observable functions of particular interest are  the Koopman eigenfunctions
\begin{math}
	\{ \varphi_1, \dots, \varphi_N \}  
\end{math}
where $\varphi_j \colon \mathcal M \to \mathbb R$ for $1\leq j \leq N$ satisfies \cite{rowleyKoop09}
\begin{equation}\label{eq:koopevp}
	\varphi_j\left({\boldsymbol{z}}_{k+1}\right) = \mathcal{K} \varphi_j\left({\boldsymbol{z}}_k\right)=\lambda_j \varphi_j\left({\boldsymbol{z}}_k\right)
\end{equation} and
\begin{math}
	\mathcal F_N := \operatorname{span} \{ \varphi_1, \dots, \varphi_N \} 
\end{math}.

To generalize  the action of the Koopman operator  vector valued observable 
\begin{equation}
	\boldsymbol{g}\colon \mathcal{M} \to \mathbb{R}^{p} \quad \boldsymbol{g}(\boldsymbol{z}) =
	\begin{bmatrix}
		g_1(\boldsymbol{z}) \\
		g_2(\boldsymbol{z}) \\
		\vdots \\
		g_p(\boldsymbol{z})
	\end{bmatrix}
\end{equation}
is considered,
whose components are given by scalar valued observable functions $g_i\colon\mathcal{M} \to \mathbb{R}$ where $g_i \in \mathcal{F}_N$. 
Using the spectral information on the Koopman operator we can expand the vector valued observable  
$\boldsymbol{g}(\boldsymbol{z}) = \sum^{N}_{j=1} \boldsymbol{v}_{j} \varphi_j(\boldsymbol{z})$,
where the vectors containing the basis coefficients $\boldsymbol{v}_j$ are referred to as the Koopman modes.

A particular choice for the vector valued observable function is the system's full state,
i.e.,
$\boldsymbol{g}(\boldsymbol{z}) = \boldsymbol{z}$.
In this case the Koopman operator perspective yields an alternative representation of the state of the system
\begin{equation}\label{eq:fullstate-eig}
	\boldsymbol{z}_k = \boldsymbol{g}(\boldsymbol{z}_k) = \sum^{N}_{j=1} \boldsymbol{v}_{j} \varphi_j(\boldsymbol{z}_k)
\end{equation}
and its temporal evolution
\begin{equation}\label{eq:fullstate-evolution-eig}
	\mathcal{K} \boldsymbol{g}(z) = \boldsymbol{g}(\boldsymbol{F}(\boldsymbol{z}_k)) = \boldsymbol{F}(\boldsymbol{z}_k) = \sum^{N}_{j=1} \boldsymbol{v}_{j} \lambda_{j} \varphi_j(\boldsymbol{z}_k),
\end{equation}
also known as Koopman mode decomposition \cite{mezicKoop05}.

\section{Sensitivity analysis}\label{sec:sens-ana-appendix}
In order to gain a better understanding of how the number of retained modes $r$ and the number of time delay embeddings $s$ affect the quality of the DMD models, a sensitivity analysis is conducted in which one of the model parameters is fixed while the other is increased. Throughout this study the ExactDMD algorithm is used and the tracked error metrics are the reconstruction error in \begin{equation}\label{eq:error-train}
	e_{\mathrm{train}} = \frac{\bigl\| \boldsymbol{Z}_{\mathrm{DMD,\,train}} - \boldsymbol{Z}_{\mathrm{train}} \bigr\|_F}{\bigl\| \boldsymbol{Z}_{\mathrm{train}} \bigr\|_F},
\end{equation}
where and the spectral agreement  
\begin{equation}\label{eq:spec-agree}
	\gamma = \frac{\langle \hat{f}^\mathrm{DMD}_i , \hat{f}_i \rangle_2}{\|\hat{f}_i\|_2\,\|\hat{f}^\mathrm{DMD}_i\|_2},
\end{equation}
where $\hat{f}_i$ and $\hat{f}^\mathrm{DMD}_i$ denote the vectors of amplitudes obtained from the FFT of probes of the original and reconstructed signal, respectively. 

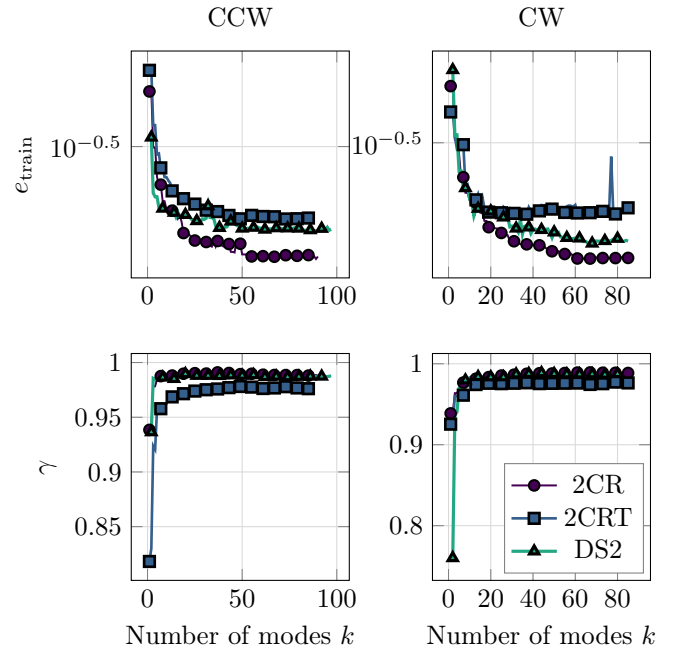
\begin{figure}[htpb]
	\centering
	\pgfplotstableread[col sep = comma]{4_7_trunkModeSens_specAgree_CCW_2CR.csv}  \tmSACCWIICR
\pgfplotstableread[col sep = comma]{4_7_trunkModeSens_specAgree_CCW_2CRT.csv} \tmSACCWIICRT
\pgfplotstableread[col sep = comma]{4_7_trunkModeSens_specAgree_CCW_DS2.csv}  \tmSACCWDS
\pgfplotstableread[col sep = comma]{4_7_trunkModeSens_specAgree_CW_2CR.csv}   \tmSACWIICR
\pgfplotstableread[col sep = comma]{4_7_trunkModeSens_specAgree_CW_2CRT.csv}  \tmSACWIICRT
\pgfplotstableread[col sep = comma]{4_7_trunkModeSens_specAgree_CW_DS2.csv}   \tmSACWDS
\pgfplotstableread[col sep = comma]{4_7_trunkModeSens_trainError_CCW_2CR.csv}  \tmETCCWIICR
\pgfplotstableread[col sep = comma]{4_7_trunkModeSens_trainError_CCW_2CRT.csv} \tmETCCWIICRT
\pgfplotstableread[col sep = comma]{4_7_trunkModeSens_trainError_CCW_DS2.csv}  \tmETCCWDS
\pgfplotstableread[col sep = comma]{4_7_trunkModeSens_trainError_CW_2CR.csv}   \tmETCWIICR
\pgfplotstableread[col sep = comma]{4_7_trunkModeSens_trainError_CW_2CRT.csv}  \tmETCWIICRT
\pgfplotstableread[col sep = comma]{4_7_trunkModeSens_trainError_CW_DS2.csv}   \tmETCWDS

\begin{tikzpicture}
   \begin{groupplot}[
        group style={
            group size=2 by 2,
            horizontal sep=11mm,
            vertical sep=10mm,
            xlabels at=edge bottom,
            ylabels at=edge left,
        },
        width=0.33\columnwidth,
        height=3cm,
        scale only axis,
        axis on top=false,
        clip=false,
        colormap/viridis,
        xlabel near ticks,
        ylabel near ticks,
        xlabel={Number of modes $k$},
        grid=both,
        grid style={gray!30},
        legend style={
            fill=white,
            draw=gray,
            inner sep=2pt
        },
    ]

    \nextgroupplot[
        ylabel={$e_\mathrm{train}$},
        ymode=log,
        title={CCW},
	cycle list name = three_col,
        legend pos=north east,
    ]
    \addplot  table[x index=0, y index=1] {\tmETCCWIICR};
    \addplot table[x index=0, y index=1] {\tmETCCWIICRT};
    \addplot  table[x index=0, y index=1] {\tmETCCWDS};

    \nextgroupplot[
        ymode=log,
        title={CW},
	cycle list name = three_col,
        legend pos=north east,
    ]
    \addplot  table[x index=0, y index=1] {\tmETCWIICR};
    \addplot table[x index=0, y index=1] {\tmETCWIICRT};
    \addplot  table[x index=0, y index=1] {\tmETCWDS};

    \nextgroupplot[
        ylabel={$\gamma$},
	cycle list name = three_col,
        legend pos=south east,
    ]
    \addplot  table[x index=0, y index=1] {\tmSACCWIICR};
    \addplot table[x index=0, y index=1] {\tmSACCWIICRT};
    \addplot  table[x index=0, y index=1] {\tmSACCWDS};

    \nextgroupplot[
        legend pos=south east,
	cycle list name = three_col,
    ]
    \addplot  table[x index=0, y index=1] {\tmSACWIICR};
    \addplot table[x index=0, y index=1] {\tmSACWIICRT};
    \addplot  table[x index=0, y index=1] {\tmSACWDS};
    \legend{2CR, 2CRT, DS2}

    \end{groupplot}
\end{tikzpicture}
	\caption{Comparison of the reconstruction error in panels (a) and (b) and the spectral agreement in panels (c) and (d) for the CCW and CW moving FR, in dependence of an increasing  number of retained modes $r$ at fixed $s = s_{\mathrm{Base}}$.}
	\label{fig:05iModesErrorMetricModes}
\end{figure}

First, the number of time delay embeddings is fixed at $s = s_{\mathrm{Base}}$ and the number of retained modes $r$ is increased from $r=1$ up to $r \in \{r_{\mathrm{CCW}},\,r_{\mathrm{CW}}\}$.
Figure~\ref{fig:05iModesErrorMetricModes}(a) and (b) show the reconstruction error for an increasing number of retained modes in the CCW and CW moving FR, respectively. 
It is observed that, as the number of retained modes increases, the reconstruction error decays in a convergence-like manner for all datasets.  
The converging decay of the error suggests that increasing the number of modes beyond the optimal values $r \in \{ r_\mathrm{CCW},\,r_\mathrm{CW}\}$ determined by the optimal threshold criterion would result in capturing sensor noise instead of the relevant part of the dynamics.

\begin{figure}[htpb]
	\centering
	\pgfplotstableread[col sep = comma]{4_7_StackedStateSens_specAgree_CCW_2CR.csv}  \saSACCWIICR 
\pgfplotstableread[col sep = comma]{4_7_StackedStateSens_specAgree_CCW_2CRT.csv} \saSACCWIICRT
\pgfplotstableread[col sep = comma]{4_7_StackedStateSens_specAgree_CCW_DS2.csv}  \saSACCWDS
\pgfplotstableread[col sep = comma]{4_7_StackedStateSens_specAgree_CW_2CR.csv}   \saSACWIICR
\pgfplotstableread[col sep = comma]{4_7_StackedStateSens_specAgree_CW_2CRT.csv}  \saSACWIICRT
\pgfplotstableread[col sep = comma]{4_7_StackedStateSens_specAgree_CW_DS2.csv}   \saSACWDS
\pgfplotstableread[col sep = comma]{4_7_StackedStateSens_trainError_CCW_2CR.csv}  \saETCCWIICR
\pgfplotstableread[col sep = comma]{4_7_StackedStateSens_trainError_CCW_2CRT.csv} \saETCCWIICRT
\pgfplotstableread[col sep = comma]{4_7_StackedStateSens_trainError_CCW_DS2.csv}  \saETCCWDS
\pgfplotstableread[col sep = comma]{4_7_StackedStateSens_trainError_CW_2CR.csv}   \saETCWIICR
\pgfplotstableread[col sep = comma]{4_7_StackedStateSens_trainError_CW_2CRT.csv}  \saETCWIICRT
\pgfplotstableread[col sep = comma]{4_7_StackedStateSens_trainError_CW_DS2.csv}   \saETCWDS

\begin{tikzpicture}
    \pgfplotsset{
        mode2CR/.style  = {color of colormap={0   of viridis}},
        mode2CRT/.style = {color of colormap={500 of viridis}},
        modeDS2/.style  = {color of colormap={900 of viridis}},
    }
    \begin{groupplot}[
        group style={
            group size=2 by 2,
            horizontal sep=11mm,
            vertical sep=10mm,
            xlabels at=edge bottom,
            ylabels at=edge left,
        },
        width=0.33\columnwidth,
        height=3cm,
        scale only axis,
        axis on top=false,
        colormap/viridis,
        xlabel near ticks,
        ylabel near ticks,
        xlabel={Number of stacks $s$},
        grid=both,
        grid style={gray!30},
        legend style={
            fill=white,
            draw=gray,
            rounded corners,
            inner sep=2pt
        },
    ]

    \nextgroupplot[
        ylabel={$e_\mathrm{train}$},
        ymode=log,
        ymax = 1,
	cycle list name = three_col,
        title={CCW},
        legend pos=north east,
    ]
    \addplot table[x index=0, y index=1] {\saETCCWIICR};
    \addplot table[x index=0, y index=1] {\saETCCWIICRT};
    \addplot  table[x index=0, y index=1] {\saETCCWDS};

    \nextgroupplot[
        ymode=log,
        ymax = 1,
        title={CW},
	cycle list name = three_col,
        legend pos=north east,
    ]
    \addplot  table[x index=0, y index=1] {\saETCWIICR};
    \addplot table[x index=0, y index=1] {\saETCWIICRT};
    \addplot  table[x index=0, y index=1] {\saETCWDS};

    \nextgroupplot[
        ylabel={$\gamma$},
	cycle list name = three_col,
        ymin=0.8, ymax=1.0,
        legend pos=south east,
    ]
    \addplot  table[x index=0, y index=1] {\saSACCWIICR};
    \addplot table[x index=0, y index=1] {\saSACCWIICRT};
    \addplot  table[x index=0, y index=1] {\saSACCWDS};

    \nextgroupplot[
        legend pos=south east,
        ymin=0.8, ymax=1.0,
	cycle list name = three_col,
    ]
    \addplot  table[x index=0, y index=1] {\saSACWIICR};
    \addplot table[x index=0, y index=1] {\saSACWIICRT};
    \addplot  table[x index=0, y index=1] {\saSACWDS};
    \legend{2CR, 2CRT, DS2}

    \end{groupplot}
\end{tikzpicture}
	\caption{Comparison of the reconstruction error in panels (a) and (b) and the spectral agreement in panels (c) and (d) for the CCW and CW moving FR, in dependence of an increasing  number of stacks $s$ at fixed $r \in \{ r_\mathrm{CCW},\,r_\mathrm{CW}\}$.}
	\label{fig:05iModesErrorMetricStacks}
\end{figure}
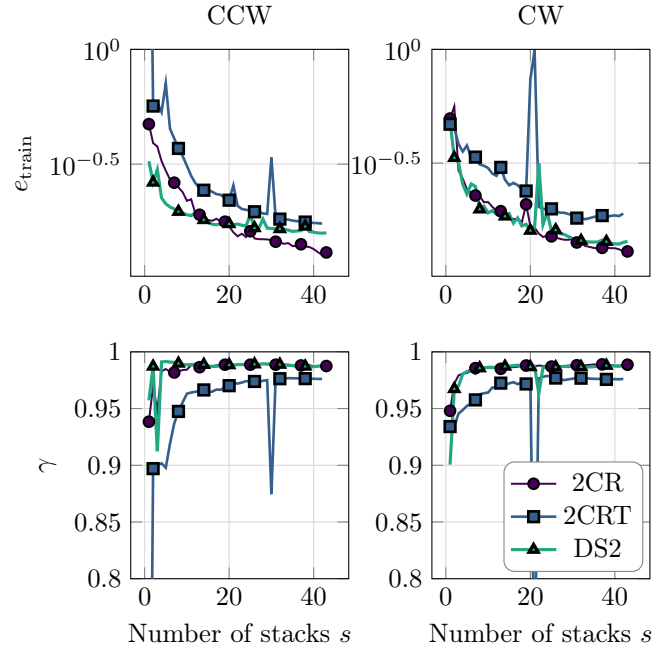

A similar behavior can be seen for the spectral agreement, which is shown in Figure~\ref{fig:05iModesErrorMetricModes}\,(c) and (d) for both moving FR.  
In particular, the high agreement is also due to the change of the FR, as one of the two counter rotating traveling waves becomes stationary in the moving FR and can therefore be described very accurately through a low number of modes.  
Subsequently, the number of modes retained is fixed at $r \in \{ r_\mathrm{CCW},\,r_\mathrm{CW}\}$, and the number of time delay embeddings of the state is increased from $s=0$ up to $s = s_\mathrm{Base}$.  
Figure~\ref{fig:05iModesErrorMetricStacks}\,(a) and (b) show the reconstruction error for an increasing number of embeddings for both moving FR.  
It is observed that, across all operating modes, increasing the number of time delay embeddings significantly improves the quality of the reconstruction.  
This improvement is partly due to the beneficial effect of time delay embedding on noise, and partly because the nonlinear phenomena can be accurately reconstructed by using the extended dictionary.  
Furthermore, the number of time delay embeddings computed from the spectral properties of the operating modes appears to be a good estimate, as the reconstruction error again decays in a convergence-like manner. This suggests that introducing additional time delay embeddings beyond this point could lead to rank deficiency, as the columns of the extended data matrices become linearly dependent.  
Finally, Figure~\ref{fig:05iModesErrorMetricStacks}\,(c) and (d) show that spectral agreement also improves with time delay embedding. However, the improvement is less pronounced, as the spectral agreement is more sensitive to the correct number of retained modes, which is already assumed.
A concluding observation is that certain combinations of the number of modes $r$ and time delay embeddings $s$ can lead to poor conditioning of the projected operator $\tilde{\boldsymbol{K}}_{\mathrm{DMD}}$.  
This results in spurious values for the reconstruction error and spectral agreement, which appear as sudden peaks in the corresponding plots.

\addcontentsline{toc}{section}{References}
\bibliographystyle{abbrvurl}
\bibliography{literature.bib}

\begin{thebibliography}{10}

\bibitem{anandRDC19}
V.~Anand and E.~Gutmark.
\newblock Rotating detonation combustors and their similarities to rocket instabilities.
\newblock {\em Progress in Energy and Combustion Science}, 73:182--234, July 2019.
\newblock \href {https://doi.org/10.1016/j.pecs.2019.04.001} {\path{doi:10.1016/j.pecs.2019.04.001}}.

\bibitem{ashkamNoiseOptDMD17}
T.~Askham and J.~N. Kutz.
\newblock Variable projection methods for an optimized dynamic mode decomposition.
\newblock {\em SIAM Journal on Applied Dynamical Systems}, 17(1):380--416, Jan. 2018.
\newblock \href {https://doi.org/10.1137/M1124176} {\path{doi:10.1137/M1124176}}.

\bibitem{bachRDC20}
E.~Bach, C.~O. Paschereit, P.~Stathopoulos, and M.~Bohon.
\newblock {{RDC}} operation and performance with varying air injector pressure loss.
\newblock In {\em {{AIAA Scitech}} 2020 {{Forum}}}, Orlando, FL, Jan. 2020. {American Institute of Aeronautics and Astronautics}.
\newblock \href {https://doi.org/10.2514/6.2020-0199} {\path{doi:10.2514/6.2020-0199}}.

\bibitem{bagheriNoiseDMD14}
S.~Bagheri.
\newblock Effects of weak noise on oscillating flows: {{Linking}} quality factor, {{Floquet}} modes, and {{Koopman}} spectrum.
\newblock {\em Physics of Fluids}, 26(9):094104, Sept. 2014.
\newblock \href {https://doi.org/10.1063/1.4895898} {\path{doi:10.1063/1.4895898}}.

\bibitem{baoCheap2RichRDE26}
Y.~Bao, J.~Zajac, M.~Powers, V.~Raman, and J.~N. Kutz.
\newblock {{Cheap2Rich}}: {{A}} multi-fidelity framework for data assimilation and system identification of multiscale physics -- rotating detonation engines, Jan. 2026.
\newblock \href {https://arxiv.org/abs/2601.20295} {\path{arXiv:2601.20295}}, \href {https://doi.org/10.48550/arXiv.2601.20295} {\path{doi:10.48550/arXiv.2601.20295}}.

\bibitem{barnouinRDC25}
P.~Barnouin, C.~O. Paschereit, E.~Bach, and M.~D. Bohon.
\newblock Dynamics and interactions of counter-rotating waves in rotating detonation combustors.
\newblock {\em Experimental Thermal and Fluid Science}, 168:111486, Sept. 2025.
\newblock \href {https://doi.org/10.1016/j.expthermflusci.2025.111486} {\path{doi:10.1016/j.expthermflusci.2025.111486}}.

\bibitem{bluemmnerRDC19}
R.~Bluemner, M.~Bohon, C.~Paschereit, and E.~Gutmark.
\newblock Counter-rotating wave mode transition dynamics in an {{RDC}}.
\newblock {\em International Journal of Hydrogen Energy}, 44(14):7628--7641, Mar. 2019.
\newblock \href {https://doi.org/10.1016/j.ijhydene.2019.01.262} {\path{doi:10.1016/j.ijhydene.2019.01.262}}.

\bibitem{bohonRDE19}
M.~Bohon, R.~Bluemner, C.~Paschereit, and E.~Gutmark.
\newblock High-speed imaging of wave modes in an {{RDC}}.
\newblock {\em Experimental Thermal and Fluid Science}, 102:28--37, Apr. 2019.
\newblock \href {https://doi.org/10.1016/j.expthermflusci.2018.10.031} {\path{doi:10.1016/j.expthermflusci.2018.10.031}}.

\bibitem{bohonDMD21}
M.~D. Bohon, A.~Orchini, R.~Bluemner, C.~O. Paschereit, and E.~J. Gutmark.
\newblock Dynamic mode decomposition analysis of rotating detonation waves.
\newblock {\em Shock Waves}, 31(7):637--649, Oct. 2021.
\newblock \href {https://doi.org/10.1007/s00193-020-00975-8} {\path{doi:10.1007/s00193-020-00975-8}}.

\bibitem{bruntonKoop21}
S.~L. Brunton, M.~Budi{\v s}i{\'c}, E.~Kaiser, and J.~N. Kutz.
\newblock Modern {{Koopman}} theory for dynamical systems.
\newblock {\em SIAM Review}, 64(2):229--340, May 2022.
\newblock \href {https://doi.org/10.1137/21M1401243} {\path{doi:10.1137/21M1401243}}.

\bibitem{bruntonBook21}
S.~L. Brunton and J.~N. Kutz.
\newblock {\em Data-Driven Science and Engineering: {{Machine}} Learning, Dynamical Systems, and Control}.
\newblock Cambridge University Press, 2 edition, May 2022.
\newblock \href {https://doi.org/10.1017/9781009089517} {\path{doi:10.1017/9781009089517}}.

\bibitem{bruntonSindy16}
S.~L. Brunton, J.~L. Proctor, and J.~N. Kutz.
\newblock Discovering governing equations from data by sparse identification of nonlinear dynamical systems.
\newblock {\em Proceedings of the National Academy of Sciences}, 113(15):3932--3937, Apr. 2016.
\newblock \href {https://doi.org/10.1073/pnas.1517384113} {\path{doi:10.1073/pnas.1517384113}}.

\bibitem{camachoProjROM2025}
R.~Camacho and C.~Huang.
\newblock Investigations on {{Projection-Based Reduced-Order Model Development}} for {{Rotating Detonation Engine}}.
\newblock {\em AIAA Journal}, 63(3):854--869, Mar. 2025.
\newblock \href {https://doi.org/10.2514/1.J064228} {\path{doi:10.2514/1.J064228}}.

\bibitem{colbrookHHDMD24}
M.~J. Colbrook.
\newblock The multiverse of dynamic mode decomposition algorithms.
\newblock In {\em Handbook of {{Numerical Analysis}}}, volume~25, pages 127--230. Elsevier, 2024.
\newblock \href {https://doi.org/10.1016/bs.hna.2024.05.004} {\path{doi:10.1016/bs.hna.2024.05.004}}.

\bibitem{dawsonNoiseDMD16}
S.~T.~M. Dawson, M.~S. Hemati, M.~O. Williams, and C.~W. Rowley.
\newblock Characterizing and correcting for the effect of sensor noise in the dynamic mode decomposition.
\newblock {\em Experiments in Fluids}, 57(3):42, Mar. 2016.
\newblock \href {https://arxiv.org/abs/1507.02264} {\path{arXiv:1507.02264}}, \href {https://doi.org/10.1007/s00348-016-2127-7} {\path{doi:10.1007/s00348-016-2127-7}}.

\bibitem{farcasOpInf2023}
I.~Farcas, R.~Gundevia, R.~Munipalli, and K.~E. Willcox.
\newblock Parametric non-intrusive reduced-order models via operator inference for large-scale rotating detonation engine simulations.
\newblock In {\em {{AIAA SCITECH}} 2023 {{Forum}}}, National Harbor, MD \& Online, Jan. 2023. {American Institute of Aeronautics and Astronautics}.
\newblock \href {https://doi.org/10.2514/6.2023-0172} {\path{doi:10.2514/6.2023-0172}}.

\bibitem{forootaniSindy24}
A.~Forootani and P.~Benner.
\newblock {{GN-SINDy}}: {{Greedy}} sampling neural network in sparse identification of nonlinear partial differential equations, May 2024.
\newblock \href {https://arxiv.org/abs/2405.08613} {\path{arXiv:2405.08613}}, \href {https://doi.org/10.48550/arXiv.2405.08613} {\path{doi:10.48550/arXiv.2405.08613}}.

\bibitem{gavishSigTresh14}
M.~Gavish and D.~L. Donoho.
\newblock The optimal hard yhreshold for singular values is \textbackslash (4/\textbackslash sqrt \textbraceleft 3\textbraceright\textbackslash ).
\newblock {\em IEEE Transactions on Information Theory}, 60(8):5040--5053, Aug. 2014.
\newblock \href {https://doi.org/10.1109/TIT.2014.2323359} {\path{doi:10.1109/TIT.2014.2323359}}.

\bibitem{goyalDeepLearn23}
P.~Goyal, S.~Y{\i}ld{\i}z, and P.~Benner.
\newblock Deep learning for structure-preserving universal stable {{Koopman-inspired}} embeddings for nonlinear canonical {{Hamiltonian}} dynamics.
\newblock {\em Machine Learning: Science and Technology}, 6(1):015063, Mar. 2025.
\newblock \href {https://doi.org/10.1088/2632-2153/adb9b5} {\path{doi:10.1088/2632-2153/adb9b5}}.

\bibitem{gutmarkPGC21}
E.~J. Gutmark.
\newblock Pressure gain combustion.
\newblock {\em Shock Waves}, 31(7):619--621, Oct. 2021.
\newblock \href {https://doi.org/10.1007/s00193-021-01053-3} {\path{doi:10.1007/s00193-021-01053-3}}.

\bibitem{hematiNoiseDMD17}
M.~S. Hemati, C.~W. Rowley, E.~A. Deem, and L.~N. Cattafesta.
\newblock De-biasing the dynamic mode decomposition for applied {{Koopman}} spectral analysis of noisy datasets.
\newblock {\em Theoretical and Computational Fluid Dynamics}, 31(4):349--368, Aug. 2017.
\newblock \href {https://doi.org/10.1007/s00162-017-0432-2} {\path{doi:10.1007/s00162-017-0432-2}}.

\bibitem{johnsonApplicationConvolutionalNeural2021}
K.~B. Johnson, D.~H. Ferguson, R.~S. Tempke, and A.~C. Nix.
\newblock Application of a {{Convolutional Neural Network}} for {{Wave Mode Identification}} in a {{Rotating Detonation Combustor Using High-Speed Imaging}}.
\newblock {\em Journal of Thermal Science and Engineering Applications}, 13(6):061021, Dec. 2021.
\newblock \href {https://doi.org/10.1115/1.4049868} {\path{doi:10.1115/1.4049868}}.

\bibitem{kailasanathPCG00}
K.~Kailasanath.
\newblock Review of propulsion applications of detonation waves.
\newblock {\em AIAA Journal}, 38(9):1698--1708, Sept. 2000.
\newblock \href {https://doi.org/10.2514/2.1156} {\path{doi:10.2514/2.1156}}.

\bibitem{kochNeuralODE2021}
J.~Koch.
\newblock Data-driven surrogates of rotating detonation engine physics with neural ordinary differential equations and high-speed camera footage.
\newblock {\em Physics of Fluids}, 33(9):091703, Sept. 2021.
\newblock \href {https://doi.org/10.1063/5.0063624} {\path{doi:10.1063/5.0063624}}.

\bibitem{kochModelRDE20}
J.~Koch, M.~Kurosaka, C.~Knowlen, and J.~N. Kutz.
\newblock Mode-locked rotating detonation waves: {{Experiments}} and a model equation.
\newblock {\em Physical Review E}, 101(1):013106, Jan. 2020.
\newblock \href {https://doi.org/10.1103/PhysRevE.101.013106} {\path{doi:10.1103/PhysRevE.101.013106}}.

\bibitem{kochMultiRDE21}
J.~Koch, M.~Kurosaka, C.~Knowlen, and J.~N. Kutz.
\newblock Multiscale physics of rotating detonation waves: {{Autosolitons}} and modulational instabilities.
\newblock {\em Physical Review E}, 104(2):024210, Aug. 2021.
\newblock \href {https://doi.org/10.1103/PhysRevE.104.024210} {\path{doi:10.1103/PhysRevE.104.024210}}.

\bibitem{koopHam31}
B.~O. Koopman.
\newblock Hamiltonian systems and transformation in {{Hilbert}} space.
\newblock {\em Proceedings of the National Academy of Sciences}, 17(5):315--318, May 1931.
\newblock \href {https://doi.org/10.1073/pnas.17.5.315} {\path{doi:10.1073/pnas.17.5.315}}.

\bibitem{luschNNDMD18}
B.~Lusch, J.~N. Kutz, and S.~L. Brunton.
\newblock Deep learning for universal linear embeddings of nonlinear dynamics.
\newblock {\em Nature Communications}, 9(1):4950, Nov. 2018.
\newblock \href {https://arxiv.org/abs/1712.09707} {\path{arXiv:1712.09707}}, \href {https://doi.org/10.1038/s41467-018-07210-0} {\path{doi:10.1038/s41467-018-07210-0}}.

\bibitem{mendibleROMwave20}
A.~Mendible, S.~L. Brunton, A.~Y. Aravkin, W.~Lowrie, and J.~N. Kutz.
\newblock Dimensionality reduction and reduced order modeling for traveling wave physics.
\newblock {\em Theoretical and Computational Fluid Dynamics}, 34(4):385--400, Aug. 2020.
\newblock \href {https://arxiv.org/abs/1911.00565} {\path{arXiv:1911.00565}}, \href {https://doi.org/10.1007/s00162-020-00529-9} {\path{doi:10.1007/s00162-020-00529-9}}.

\bibitem{mendibleRDE21}
A.~Mendible, J.~Koch, H.~Lange, S.~L. Brunton, and J.~N. Kutz.
\newblock Data-driven modeling of rotating detonation waves.
\newblock {\em Physical Review Fluids}, 6(5):050507, May 2021.
\newblock \href {https://doi.org/10.1103/PhysRevFluids.6.050507} {\path{doi:10.1103/PhysRevFluids.6.050507}}.

\bibitem{mezicKoop05}
I.~Mezi{\'c}.
\newblock Spectral properties of dynamical systems, model deduction and decompositions.
\newblock {\em Nonlinear Dynamics}, 41(1-3):309--325, Aug. 2005.
\newblock \href {https://doi.org/10.1007/s11071-005-2824-x} {\path{doi:10.1007/s11071-005-2824-x}}.

\bibitem{williamsKerDMD15}
M.~O.~Williams, C.~W.~Rowley, and I.~G.~Kevrekidis.
\newblock A kernel-based method for data-driven {{Koopman}} spectral analysis.
\newblock {\em Journal of Computational Dynamics}, 2(2):247--265, 2015.
\newblock \href {https://doi.org/10.3934/jcd.2015005} {\path{doi:10.3934/jcd.2015005}}.

\bibitem{oexleHSV2026}
D.~Oexle and M.~Bohon.
\newblock High speed video image and pressure data of counter rotating detonation waves, Mar. 2026.
\newblock \href {https://doi.org/10.5281/ZENODO.18886925} {\path{doi:10.5281/ZENODO.18886925}}.

\bibitem{rustVarPro13}
D.~P. O'Leary and B.~W. Rust.
\newblock Variable projection for nonlinear least squares problems.
\newblock {\em Computational Optimization and Applications}, 54(3):579--593, Apr. 2013.
\newblock \href {https://doi.org/10.1007/s10589-012-9492-9} {\path{doi:10.1007/s10589-012-9492-9}}.

\bibitem{ramanNonRDE2023}
V.~Raman, S.~Prakash, and M.~Gamba.
\newblock Nonidealities in rotating detonation engines.
\newblock {\em Annual Review of Fluid Mechanics}, 55(1):639--674, Jan. 2023.
\newblock \href {https://doi.org/10.1146/annurev-fluid-120720-032612} {\path{doi:10.1146/annurev-fluid-120720-032612}}.

\bibitem{reissSPOD18}
J.~Reiss, P.~Schulze, J.~Sesterhenn, and V.~Mehrmann.
\newblock The shifted proper orthogonal decomposition: {{A}} mode decomposition for multiple transport phenomena.
\newblock {\em SIAM Journal on Scientific Computing}, 40(3):A1322--A1344, Jan. 2018.
\newblock \href {https://doi.org/10.1137/17M1140571} {\path{doi:10.1137/17M1140571}}.

\bibitem{rowleyKoop09}
C.~W. Rowley, I.~Mezi{\'c}, S.~Bagheri, P.~Schlatter, and D.~S. Henningson.
\newblock Spectral analysis of nonlinear flows.
\newblock {\em Journal of Fluid Mechanics}, 641:115--127, Dec. 2009.
\newblock \href {https://doi.org/10.1017/S0022112009992059} {\path{doi:10.1017/S0022112009992059}}.

\bibitem{sashidarBOPDMD22}
D.~Sashidhar and J.~N. Kutz.
\newblock Bagging, optimized dynamic mode decomposition ({{BOP-DMD}}) for robust, stable forecasting with spatial and temporal uncertainty-quantification.
\newblock {\em Philosophical Transactions of the Royal Society A: Mathematical, Physical and Engineering Sciences}, 380(2229):20210199, Aug. 2022.
\newblock \href {https://arxiv.org/abs/2107.10878} {\path{arXiv:2107.10878}}, \href {https://doi.org/10.1098/rsta.2021.0199} {\path{doi:10.1098/rsta.2021.0199}}.

\bibitem{schmidDMD10}
P.~J. Schmid.
\newblock Dynamic mode decomposition of numerical and experimental data.
\newblock {\em Journal of Fluid Mechanics}, 656:5--28, Aug. 2010.
\newblock \href {https://doi.org/10.1017/S0022112010001217} {\path{doi:10.1017/S0022112010001217}}.

\bibitem{shankRDC12}
J.~C. Shank.
\newblock Development and testing of a rotating detonation engine run on hydrogen and air.
\newblock Master's thesis, Air Force Institute of Technology, Mar. 2012.

\bibitem{tuDMD14}
J.~H. Tu, C.~W. Rowley, D.~M. Luchtenburg, S.~L. Brunton, and J.~N. Kutz.
\newblock On dynamic mode decomposition: {{Theory}} and applications.
\newblock {\em Journal of Computational Dynamics}, 1(2):391--421, 2014.
\newblock \href {https://arxiv.org/abs/1312.0041} {\path{arXiv:1312.0041}}, \href {https://doi.org/10.3934/jcd.2014.1.391} {\path{doi:10.3934/jcd.2014.1.391}}.

\bibitem{williamsEDMD15}
M.~O. Williams, I.~G. Kevrekidis, and C.~W. Rowley.
\newblock A data--driven approximation of the {{Koopman}} operator: Extending dynamic mode decomposition.
\newblock {\em Journal of Nonlinear Science}, 25(6):1307--1346, Dec. 2015.
\newblock \href {https://doi.org/10.1007/s00332-015-9258-5} {\path{doi:10.1007/s00332-015-9258-5}}.

\end{thebibliography}
\end{document}